\UseRawInputEncoding 
\documentclass[a4paper]{article}
\usepackage[margin=1in]{geometry}
\usepackage{latexsym,enumerate}
\usepackage{eepic}
\usepackage{epic}
\usepackage{graphicx}
\usepackage{color}
\usepackage{ifpdf}
\usepackage{amssymb,amsmath,epsfig,algorithm,algorithmicx,multirow}
\usepackage{amsfonts,dsfont}
\usepackage{subfigure}
\usepackage{multirow}
\usepackage{makecell}
\usepackage{bm}
\usepackage{tikz}
\usepackage{xcolor}
\usetikzlibrary{graphs, positioning, quotes, shapes.geometric}

\newtheorem{remark}{Remark}[section]

\title{\Large \bf Stochastic Domain Decomposition Based on Variable-Separation Method
\thanks{The authors are grateful to the anonymous referees and the handling editor for their helpful suggestions.
	The research of this work was supported by the National Key R\&D Program of China (No. 2021YFA1001300),
	the National Natural Science Foundation of China (Nos. 12271150, 12101216, 12171406),  the Natural Science Foundation of Hunan Province (Nos. 2023JJ10001, 2022RC1190,  2022JJ40030),
	and the Hong Kong RGC grant projects (Nos. 17300318, 17307921).}
}

\author{
Liang Chen\thanks{School of Mathematical, Hunan University, Changsha 410082, China  ({\tt chl@hnu.edu.cn}).}
\and
Yaru Chen\thanks{School of Mathematical, Hunan University, Changsha 410082, China  ({\tt yrchen@hnu.edu.cn}).}
\and
Qiuqi Li\thanks{School of Mathematical, Hunan University, Changsha 410082, China  ({\tt qli28@hnu.edu.cn}).}
\and
Zhiwen Zhang \thanks{Department of Mathematics, The University of Hong Kong, Pokfulam Road, Hong Kong Special Administrative Region of China ({\tt zhangzw@hku.hk}).}
}

\begin{document}
\maketitle

\begin{abstract}
This work proposes a stochastic domain decomposition method for solving steady-state partial differential equations (PDEs) with random inputs.
Specifically, based on the efficiency of the variable-separation (VS) method in simulating stochastic partial differential equations (SPDEs), we extend it to stochastic algebraic systems and employ it in the context of stochastic domain decomposition.
The resulting method, termed stochastic domain decomposition based on variable-separation (SDD-VS) alleviates the challenge commonly known as the ``curse of dimensionality" notably by leveraging explicit representations of stochastic functions derived from physical systems.
The primary objective of the proposed SDD-VS method is to obtain a separate representation of the solution for the stochastic interface problem.
To enhance computational efficiency, we introduce a two-phase approach consisting of offline and online computation.
In the offline phase, we establish an affine representation of stochastic algebraic systems by systematically applying the VS method.
During the online phase, we estimate the interface unknowns of SPDEs using a quasi-optimal separated representation, facilitating the construction of efficient surrogate models for subproblems.
We substantiate the effectiveness of our proposed approach through numerical experiments involving three specific instances, demonstrating its capability to provide accurate solutions.

\medskip
\noindent
\textbf{keywords:}
stochastic partial differential equation; stochastic domain decomposition; stochastic interface problem; variable-separation method; uncertainty quantification

\end{abstract}

\section{Introduction}
In the realm of computational science and engineering, many problems such as modeling water flow and solute transport in heterogeneous soil and aquifer formations involve uncertainties due to inadequate knowledge about physical properties and measurement noise. To provide accurate and reliable predictions, these uncertainties are often represented by random variables, and their impact on the system is explored through the lens of stochastic partial differential equations (SPDEs).
A few numerical methods have been proposed for simulating SPDEs, including Galerkin projections \cite{babuvska2005solving,frauenfelder2005finite, matthies2005galerkin}, stochastic interpolation \cite{babuvska2007stochastic, nobile2008sparse, xiu2007efficient, xiu2005high} and the methods based on deep learning \cite{ZHU2018415,wang2020a}.
This work focuses on the domain decomposition method (DDM), a promising approach for addressing stochastic problems. We aim to apply DDM techniques for solving SPDEs efficiently and enabling robust predictions in the presence of uncertainties.

Domain decomposition methods have a long history of successful application in solving deterministic problems \cite{Papadrakakis2011a, T1994ddm, IM1993ddp, B1996ddm, HO1996ddm}. These methods are to partition a large computational domain into several subdomains so that the subproblems on each subdomain can be solved independently \cite{AA1999ddm, Toselli2005domain}. Under this setting, { DDM has primarily focused on developing parallel solvers for deterministic partial differential equations (PDEs) to improve computational efficiency\cite{Cai1993an, M1997on, sun2010parallel}}. These methods are generally categorized into overlapping (Schwarz iteration) methods and non-overlapping (Schur complement) methods \cite{Chan1992on}. Our proposed stochastic domain decomposition (SDD) method is based on the latter, specifically the Schur complement methods \cite{Chan1987analysis, Mandel1993BalancingDD, Matthias2006neumann, J2001efficient}.
In this approach, the computational domain is divided into non-overlapping subdomains, and the degrees of freedom on each subdomain are separated into interior and interface parts. The global system can be reduced to a Schur complement system by applying block Gaussian elimination. The unknown variables on the interface are then determined by solving the Schur complement system.

The application of domain decomposition methods in stochastic simulations has emerged as a promising research area. Notably, Sarkar et al. \cite{Sarker2009dd} proposed a domain decomposition approach with Schur complement in the physical space and functional decomposition in the probability space for solving stochastic partial differential equations.
Subsequently, Subber et al. \cite{waad2014pre,waad2014dd} extended this method by incorporating a typical preconditioner to determine the coefficients of the interface solution.
In a separate study by Chen et al. \cite{Chen2015local}, the authors employed a domain decomposition approach to solve stochastic elliptic PDEs. They approximated the local solution in each subdomain within a low-dimensional parametric space.
In Hadigol et al. \cite{Mo2014inter}, a stochastic model reduction approach based on low-rank separated representations for the stochastic space was studied, where the coefficients of the interface solution are done by finite element tearing and interconnecting method.
Liao and Willcox \cite{Liao2015ddua} proposed an offline-online approach that combines DDM with importance sampling. This enables the use of different strategies in local systems, and the interface solution is generated by weighting precomputed PDE solutions.
{To handle} elliptic PDEs in high-contrast random media, Hou et al. \cite{Hou2017exploring} incorporated multiscale finite element methods into the DDM framework.
Incorporating a Schwarz-type iterative algorithm, Zhang et al. \cite{Zhang2018SDD}  introduced a moment-minimizing interface condition to match the stochastic interface solution.
Lastly,  Mu and Zhang \cite{Mu2018ddrd} recently combined model reduction methods with sparse polynomial approximation. They developed a stochastic domain decomposition method comprising an offline procedure and an online procedure for linear steady-state convection-diffusion equations with random coefficients.

In this paper, we focus on integrating Schur complement methods with the variable-separation (VS) method for stochastic simulations. The VS method is one of the most effective model order reduction methods and has been successfully used to solve SPDEs in a low-dimensional manifold \cite{ Li2017a, jiang2018model, Li2020a}. Building upon the extension of the VS method to stochastic algebraic systems, we propose a novel model reduction approach for solving the global stochastic interface problem of SPDEs, termed the stochastic domain decomposition based on the variable-separation method (SDD-VS).
To enhance simulation efficiency, an offline-online computational decomposition is employed for the stochastic interface problem. The offline phase consists of three stages. Firstly, the original domain is partitioned into non-overlapping subdomains, and the Schur complement system is established, including local Schur complement matrices and corresponding right-hand side vectors with random inputs. Secondly, a reduced stochastic algebraic system is generated using low-rank representations of the global stochastic Schur complement matrices and right-hand side vectors. Although the reduced model requires less computational effort compared to the original stochastic Schur complement system, it still involves the discrete degrees of the full model and may not be considered small-scale.
To address this, the third stage introduces a functional decomposition expression for the stochastic solution on the interface, achieved through the extended VS method. In the online phase, each realization of the stochastic interface problem is recovered using the outputs from the offline phase. Finally, efficient surrogate models for the stochastic subproblems can be obtained using the VS method or other model reduction techniques \cite{Tamellini2014Model,zhang2015a,li2020a}.

The SDD-VS method proposed in this paper shares the advantages of both the variable separation and the domain decomposition method.
Specifically, the SDD-VS method can alleviate the ``curse of dimensionality'' in an effective way and this is achieved by the explicit representation of the stochastic functions deduced from the physical system instead of using a suitable set of basis functions (e.g., polynomial chaos basis and radial basis functions) of the random variables.
Meanwhile, the proposed method is much easier than other stochastic domain decomposition methods when implemented, thanks to the applications of the VS method.  Furthermore, the whole computation of the SDD-VS method consists of an offline stage and an online stage. The online stage has high computational efficiency and its computational cost is completely independent of spatial discretization.
Moreover, the SDD-VS method can reduce the dimension of the random variable for local problems in subdomains, it leads to a more efficient surrogate model for each subproblem than that obtained on the entire domain.  Finally and most importantly, the proposed method maintains the same merits as DDM for deterministic PDEs, such as the ability to solve subproblems independently and in parallel.
In summary, the SDD-VS method offers an efficient and reliable approximation for SPDEs, which is particularly valuable in many-query contexts such as optimization, control design, and inverse analysis.

This paper is organized as follows. Section \ref{ssec:prelim} provides the necessary notation and preliminaries, and introduces the VS method for stochastic algebraic systems. It also provides a brief overview of the domain decomposition method for deterministic PDEs.
In Section \ref{sec-sddm}, we present the SDD-VS method proposed in this paper, detailing its key steps and procedures, and we discuss the corresponding numerical methods of the subproblems in Section \ref{sec-subproblems}.
In Section \ref{sec-numerical examples}, three numerical examples are presented to demonstrate the performance and computational advantages of the proposed SDD-VS method. Finally, we draw conclusions and offer some final remarks on the method and its potential applications.

\section{Preliminaries and notations}
\label{ssec:prelim}
In this section, we present some preliminaries and specify the notation of this paper.
Let $(\Omega,\mathcal{B},P)$ be a finite dimensional probability space, where $\Omega$ is the event space, $\mathcal{B}$ is a $\sigma$-algebra on $\Omega$, and $P$ is the probability measure on $\mathcal{B}$.
Let $D$ be a given convex and bounded physical domain with Lipschitz continuous boundary $\partial D$.
We use $\mathcal{V}$ to denote a Hilbert space defined on $D$ with an inner product defined by $(\cdot,\cdot)_{\mathcal{V}}$, and the induced norm is given by $\|\cdot\|_{\mathcal{V}}=\sqrt{(\cdot,\cdot)_{\mathcal{V}}}$.
We consider the following SPDE defined on $D$
\begin{eqnarray}
	\label{eq-SPDE}
	\left\{
	\begin{aligned}
		\mathcal{L}(x,\bm{\xi};u(x,\bm{\xi}))&=f(x,\bm{\xi}), \ \forall ~x \in D, ~ \bm{\xi} \in \Omega,\\	
		\mathcal{B}(x,\bm{\xi};u(x,\bm{\xi}))&=g(x,\bm{\xi}), \ \forall ~x \in \partial D, ~ \bm{\xi} \in \Omega,
	\end{aligned}
	\right.
\end{eqnarray}
where $\bm{\xi}\in\Omega$ is a set of real-valued random variables, $\mathcal{L}$ is a stochastic differential operator, $\mathcal{B}$ is the boundary condition operator, $f$ is the source team and $g$ is the boundary team, $u(x,\bm{\xi})$ is the solution to this stochastic PDE.
A practical instance of \eqref{eq-SPDE} is a stochastic diffusion equation given by
\begin{eqnarray*}
	\left\{
	\begin{aligned}
		-\nabla\cdot (c(x;\bm{\xi}) \nabla u(x;\bm{\xi}))&=f(x;\bm{\xi}),  ~\forall ~x \in D, ~ \bm{\xi} \in \Omega,\\	
		u(x;\bm{\xi})&=0,  ~\forall ~x \in \partial D, ~ \bm{\xi} \in \Omega,
	\end{aligned}
	\right.
\end{eqnarray*}
where $c(x;\bm{\xi})$ is the diffusion coefficient.

The weak formulation of the problem (\ref{eq-SPDE}) reads as follows: find  $u\in \mathcal{V}$ such that
\begin{eqnarray}
	\label{eq-SPDE-weak}
	a(u(\bm{\xi}),v;\bm{\xi})=b(v;\bm{\xi}),\ \forall ~v\in \mathcal{V},
\end{eqnarray}
where $a(\cdot,\cdot;\cdot)$ and $b(\cdot;\cdot)$ are a bilinear form and a linear form on $\mathcal{V}$, respectively. We assume that $a(\cdot,\cdot;\bm{\xi})$ and $b(\cdot;\bm{\xi})$ are affine with respect to $\bm{\xi}$, i.e.,
\begin{eqnarray}
	\label{eq-vs-form}
	\left\{
	\begin{aligned}
		a(w,v;\bm{\xi})&=\sum_{k=1}^{m_a}p^{k}(\bm{\xi})a^{k}(w,v), \ \forall ~w,v\in\mathcal{V},  \forall ~\bm{\xi}\in\Omega,\\
		b(v;\bm{\xi})&=\sum_{k=1}^{m_b}q^{k}(\bm{\xi})b^{k}(v),\ \forall ~v\in\mathcal{V}, \forall ~\bm{\xi}\in\Omega,
	\end{aligned}
	\right.
\end{eqnarray}
where $p^{k}(\bm{\xi}):\Omega\to R$ is a $\bm{\xi}$-dependent stochastic function and  $a^{k}:\mathcal{V}\times\mathcal{V}\to R$ is a bilinear form independent of $\bm{\xi}$, for each $k=1,\cdots,m_a$. Each $q^{k}(\bm{\xi}):\Omega\to R$ is a $\bm{\xi}$-dependent stochastic function and each $b^{k}:\mathcal{V}\to R$ is a linear form independent of $\bm{\xi}$, for $k=1,\cdots,m_b$. When $a(\cdot,\cdot;\bm{\xi})$ and $b(\cdot;\bm{\xi})$ are not affine with respect to $\bm{\xi}$, we can use the novel VS method for multivariate function \cite{Li2017a} to obtain such an affine expansion approximation for them.

In particular, we consider the finite element (FE) approximation of problem (\ref{eq-SPDE}) in an $n$-dimensional subspace $\mathcal{V}_h\subset \mathcal{V}$. Let $\{\psi_{i}\}_{i=1}^{n}$ be the set of basis functions of the FE space $\mathcal{V}_h$, the solution $u_h(\bm{\xi})$ can be represented by
$$
u_h(\bm{\xi})=\sum_{i=1}^{n}u_{i}(\bm{\xi}) \psi_{i}.
$$
With the assumption (\ref{eq-vs-form}) of affine decomposition, we have the matrix form of equation (\ref{eq-SPDE-weak}) in an FE space $\mathcal{V}_h$ as follows
\begin{equation}
	\label{eq-matrixequations}
	\Big(\sum_{k=1}^{m_{a}}p^{k}(\bm{\xi})A^{k}\Big)\mathbf{u}(\bm{\xi})=\sum_{k=1}^{m_b}q^{k}(\bm{\xi})F^{k},
\end{equation}
where
\begin{eqnarray*}
	(A^{k})_{ij}=a^{k}(\psi_{i},\psi_{j}),\quad (\mathbf{u}(\bm{\xi}))_{j}=u_{j}(\bm{\xi}), \quad (F^{k})_{j}=b^{k}(\psi_{j}),~~1\leq i,j\leq n.
\end{eqnarray*}

Now we consider equation (\ref{eq-matrixequations}) as an example and introduce the VS method for stochastic algebraic systems. This is the key technology for the stochastic domain decomposition method we will propose in Section \ref{sec-sddm}.

\subsection{VS method for stochastic algebraic systems}
\label{sec-vs}
We attempt to achieve an approximation of the stochastic algebraic system (\ref{eq-matrixequations}) in the form
\begin{eqnarray}
	\label{eq-approx}
	\mathbf{u}(\bm{\xi})\approx \mathbf{u}_{N}(\bm{\xi}):=\sum_{i=1}^{N}\zeta_i(\bm{\xi}) \mathbf{c}_i,
\end{eqnarray}
where $\zeta_i(\bm{\xi})$ are stochastic functions and $\mathbf{c}_i$ are deterministic vectors, and $N$ is the number of the separated terms. The VS method was first proposed in \cite{Li2017a} for linear stochastic problems, and has seen extensions to stochastic saddle point problems in \cite{jiang2018model},  and to nonlinear parameterized PDEs in \cite{Li2020a}. In this contribution, we develop the strategy of VS for stochastic algebraic systems. The VS method employs an offline-online computational decomposition to enhance efficiency. In offline stage, to generate the reduced basis functions $\{\mathbf{c}_i\}_{i=1}^N$ and $\{\zeta_i(\bm{\xi})\}_{i=1}^N$,  we need to compute a set of snapshots, which are the solutions of the stochastic problem (\ref{eq-matrixequations}) corresponding to a set of optimal parameter samples. In the online stage, the output is computed by the quasi-optimal separated representations for many instances of parameters, and the influence of the uncertainty is estimated.

Here we describe the detail of the VS algorithm to obtain $\{\mathbf{c}_i\}_{i=1}^N$ and $\{\zeta_i(\bm{\xi})\}_{i=1}^N$ in (\ref{eq-approx}). To this end, we define the residual for the VS method by
\begin{eqnarray}
	\label{eq-resi}
	{e}(\bm{\xi}):=\mathbf{u}(\bm{\xi})-\mathbf{u}_{k-1}(\bm{\xi}).
\end{eqnarray}
By (\ref{eq-matrixequations}), we have that
\begin{eqnarray}
	\label{resi-eq}
	\Big(\sum_{j=1}^{m_a}p^{j}(\bm{\xi})A^{j}\Big){e}(\bm{\xi})=\sum_{j=1}^{m_b}q^{j}(\bm{\xi})F^{j}-\Big(\sum_{j=1}^{m_{a}}p^{j}(\bm{\xi})A^{j}\Big)\mathbf{u}_{k-1}(\bm{\xi}).	
\end{eqnarray}
Let ${r}_k(\bm{\xi})$ be the residual of equation (\ref{eq-matrixequations}) when using $\mathbf{u}_{k-1}(\bm{\xi})$ to approximate $\mathbf{u}(\bm{\xi})$, that is,
\begin{eqnarray}
	\label{resi-eq2}
	{r}_k(\bm{\xi}):=
	\left\{
	\begin{aligned}
		&\sum_{j=1}^{m_b}q^{j}(\bm{\xi})F^{j},&k=1,\\
		&\sum_{j=1}^{m_b}q^{j}(\bm{\xi})F^{j}-\Big(\sum_{j=1}^{m_{a}}p^{j}(\bm{\xi})A^{j}\Big)\mathbf{u}_{k-1}(\bm{\xi}),&k\geq2.
	\end{aligned}
	\right.
\end{eqnarray}
Combining the representation of residual equation (\ref{resi-eq2}) with equation (\ref{resi-eq}), we get the following error residual equation as
\begin{eqnarray}
	\label{resi-eq3}
	\sum_{j=1}^{m_a}p^{j}(\bm{\xi})A^{j}{e}(\bm{\xi})={r}_k(\bm{\xi}).
\end{eqnarray}
At step $k$, we choose $\bm{\xi}_k$ as follows
\begin{eqnarray*}
	\bm{\xi}_k:=
	\left\{
	\begin{aligned}
		&\text{chosen randomly in } \Omega,&k=1,\\
		&\text{argmax}_{\bm{\xi}\in\Xi}\|{r}_k(\bm{\xi})\|_2,&k\geq2,
	\end{aligned}
	\right.
\end{eqnarray*}
where $\Xi$ is a collection of a finite number of samples in $\Omega$.
Let $\mathbf{e}_h$ be the solution of (\ref{resi-eq3}) with $\bm{\xi}=\bm{\xi}_k$, then we obtain the { $k$-th deterministic column vector $\mathbf{c}_k=\mathbf{e}_h$ in (\ref{eq-approx})}. Given ${e}(\bm{\xi}):=\mathbf{c}_ke_{\bm{\xi}}(\bm{\xi})$, and $\mathbf{u}_{k-1}(\bm{\xi}):=\sum_{i=1}^{k-1}\zeta_i(\bm{\xi})\mathbf{c}_i$ in equation (\ref{resi-eq}), it follows that
\begin{eqnarray}
	\label{resi-eq1}
	\sum_{j=1}^{m_a}p^{j}(\bm{\xi})A^{j}\mathbf{c}_k{e}_{\bm{\xi}}(\bm{\xi})=\sum_{j=1}^{m_b}q^{j}(\bm{\xi})F^{j}-\sum_{j=1}^{m_a}\sum_{i=1}^{k-1}p^{j}(\bm{\xi})\zeta_i(\bm{\xi})A^{j}\mathbf{c}_i.	
\end{eqnarray}
Both sides of the equation (\ref{resi-eq1}) are taken dot product with $\mathbf{c}_k$. Then, we have
\begin{eqnarray*}
	\label{separation spde}
	\sum_{j=1}^{m_a}p^{j}(\bm{\xi})(\mathcal{A})_{kj}{e}_{\bm{\xi}}(\bm{\xi})=\sum_{j=1}^{m_b}q^{j}(\bm{\xi})(\mathcal{F})_j-\sum_{j=1}^{m_a}\sum_{i=1}^{k-1}p^{j}(\bm{\xi})\zeta_i(\bm{\xi})(\mathcal{A})_{ij},
\end{eqnarray*}
where the matrix $\mathcal{A}$ is defined by $(\mathcal{A})_{ij}=(\mathbf{c}_k)^{T}A^{j}\mathbf{c}_i$, for $1\leq i\leq k$, $1\leq j\leq m_a$, and the vector $\mathcal{F}$ is defined by $(\mathcal{F})_j=(\mathbf{c}_k)^{T}F^{j}$, for $1\leq j\leq m_b$.
This gives rise to
\begin{eqnarray}
	\label{eq-stochatic}
	{e}_{\bm{\xi}}(\bm{\xi})=\frac{\sum_{j=1}^{m_b}q^{j}(\bm{\xi})(\mathcal{F})_j-\sum_{j=1}^{m_a}\sum_{i=1}^{k-1}p^{j}(\bm{\xi})\zeta_i(\bm{\xi})(\mathcal{A})_{ij}}{\sum_{j=1}^{m_a}p^{j}(\bm{\xi})(\mathcal{A})_{kj}}.
\end{eqnarray}
Consequently, the $k$-th stochastic function in (\ref{eq-approx}) can be obtained by taking  $\zeta_k(\bm{\xi})=e_{\bm{\xi}}(\bm{\xi})$.
When $\left\|r_k(\bm{\xi}_k)\right\|_2:=\sqrt{(r_k(\bm{\xi}_k))^{T}r_k(\bm{\xi}_k)}$ is small enough, we can stop the iteration procedure.
Algorithm \ref{algorithm-vs-sSAS} summarizes the above procedure of the VS method for stochastic algebraic systems.
\begin{algorithm}
	\caption{The VS method for stochastic algebraic systems}
	\textbf{Input}:  The stochastic algebraic system in (\ref{eq-matrixequations}), a set of samples $\Xi\in\Omega$,\\
	and the error tolerance $\varepsilon$.\\
	\textbf{Output}: The separated representation $\mathbf{u}_N(\bm{\xi}):= \sum_{i=1}^N\zeta_i(\bm{\xi})\mathbf{c}_i$.\\
	~1:~~Initialize the residual $r(\bm{\xi}):=\sum_{j=1}^{m_b}q^{j}(\bm{\xi})F^{j}$, a random $\bm{\bm{\xi}_1}\in \Xi$,\\
	$~~~~~$the iteration counter $k=1$;\\
	~2:~~Calculate $\mathbf{c}_k=\mathbf{e}_h$ by solving (\ref{resi-eq3}) with $\bm{\xi}=\bm{\bm{\xi}_k}$, $\zeta_k(\bm{\xi})=e_{\bm{\xi}}(\bm{\xi})$ by (\ref{eq-stochatic});\\
	~3:~~Update $\Xi=\Xi \textbf{\textbackslash} \bm{\xi}_k $, and take the approximation $\mathbf{u}_k(\bm{\xi}):= \sum_{i=1}^k \zeta_i(\bm{\xi})\mathbf{c}_i$;\\
	~4:~~$k\to k+1$;\\
	~5:~~Take the residual $r(\bm{\xi}):=\sum_{j=1}^{m_b}q^{j}(\bm{\xi})F^{j}-\sum_{j=1}^{m_a}p^{j}(\bm{\xi})A^{j}\mathbf{u}_{k-1}(\bm{\xi})$, and choose\\ $~~~~~$ $\bm{\xi}_k=\text{argmax}_{\bm{\xi}\in\Xi}\|r(\bm{\xi})\|_{2}$;\\
	~6:~~Return to Step 2 if $\|r(\bm{\xi})\|_{2} \geq \varepsilon$, otherwise \textbf{terminate}.\\
	~7:~~$N=k$.
	\label{algorithm-vs-sSAS}
\end{algorithm}

Note that the matrix $\mathcal{A}$ and the vector $\mathcal{F}$ are independent of random variables $\bm{\xi}$, and their computation is once in the offline phase. The online computation is to calculate (\ref{eq-approx}) for any $\bm{\xi} \in \Omega$, which is efficient because the online computation only involves the separated representation (\ref{eq-approx}).

\subsection{Domain decomposition method for deterministic problem}
\label{sec-ddm}
In this section, we briefly review the non-overlapping domain decomposition method for deterministic PDEs (see, e.g., \cite{waad2014dd, Mu2018ddrd} for more details), and take the equation (\ref{eq-SPDE}) with Dirichlet boundary condition for a realization of random variable $\bm{\xi}\in\Omega$ as an example. In this case, we have the following weak formulation which is independent of random variables
\begin{eqnarray*}
	a(u,v)=b(v),\ \forall ~v\in \mathcal{V}.
\end{eqnarray*}

\begin{figure}[htbp]
	\centering
	\begin{tikzpicture}[scale=1]
		\draw[] (0,0) ellipse (2.5 and 1.5);
		\draw [](0, 1.5) .. controls (1,1) and (-1.5,-0.5) .. (0,-1.5);
		\node[] (B) at (1,0) {$D_2$};
		\node[] (B) at (-1.2,0) {$D_1$};
		\node[] (B) at (0,-0.2) {$\Gamma_{12}$};
	\end{tikzpicture}
	\caption{Illustration of a domain partitioned into two subdomains.}
	\label{ddm}
\end{figure}
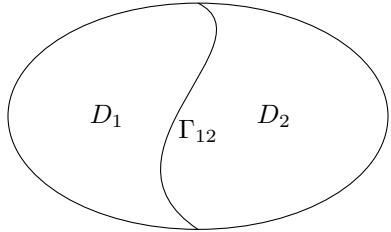

Assumed that $D$ is partitioned into $N_s$ non-overlapping subdomains, denoted by $D_i,i=1,2,\cdots,N_s$, such that
\begin{eqnarray*}
	D=\mathop{\large\cup}\limits_{i=1}^{N_s}\overline{D_i} ~~\text{and}~~ D_i\cap D_j=\emptyset, ~~\text{if}~ i\neq j,
\end{eqnarray*}
and we define the interface of two adjoint subdomains $D_i$ and $D_j$ by $\Gamma_{ij}$, a notional example for two subdomains is shown in Figure \ref{ddm}. Moreover, we denote the Hilbert space $\mathcal{V}$ restricted on subdomain $D_i$ by $\mathcal{V}_{i}$.

Then the weak formulation in a typical subdomain $D_i$ can be expressed as $a(u,v)_i=b(v)_i,\ \forall ~v\in \mathcal{V}_i$. The finite element approximation of the above equation leads to a local linear system as follows
\begin{eqnarray}
	\label{eq-local-matrix}
	A_i\bm{u}_i=\bm{f}_i,
\end{eqnarray}
where $A_i$ and $\bm{f}_i$ are the local stiffness matrix and local load vector, respectively, and $\bm{u}_i$ is the vector of local nodal values in the subdomain $D_i$.
Due to the lack of boundary conditions, the local system (\ref{eq-local-matrix}) is singular.
Then we consider dividing $\bm{u}_i$ into two parts: the nodal shared by two or more adjacent subdomains, i.e. the interface part $\bm{u}_{\Gamma}^i$ and the interior part $\bm{u}_{I}^i$ that belong to subdomain $D_i$. Consequently, the system (\ref{eq-local-matrix}) can be rewritten as
\begin{eqnarray*}
	\begin{bmatrix}A_{II}^i& A_{I\Gamma}^i \\ \\ A_{\Gamma I}^i&  A_{\Gamma\Gamma}^i
	\end{bmatrix}
	\begin{Bmatrix} \bm{u}_{I}^i\\\quad \\ \bm{u}_{\Gamma}^i  \end{Bmatrix}
	=
	\begin{Bmatrix} \bm{f}_{I}^i\\ \quad \\ \bm{f}_{\Gamma}^i \end{Bmatrix}.
\end{eqnarray*}

Performing the Gaussian elimination technique, we know that once the interface unknowns $\bm{u}_{\Gamma}^i$ are obtained, the interior unknowns $\bm{u}_{I}^i$ can be obtained by solving the interior problem on each subdomain $D_i$ as follows
\begin{eqnarray}
	\label{eq-interior}
	\bm{u}_{I}^i=(A_{II}^i)^{-1}(\bm{f}_{I}^i-A_{I\Gamma}^i\bm{u}_{\Gamma}^i).
\end{eqnarray}

Next, we discuss how to get the interface unknowns $\bm{u}_{\Gamma}^i$. First, we introduce the restriction matrix $R_i$ { (consisting of zeros and ones) represents a scatter operator}, which relates the global interface unknowns $\bm{u}_{\Gamma}$ and the local interface unknowns $\bm{u}_{\Gamma}^i$ as
\begin{eqnarray*}
	\label{global-local}
	R_i\bm{u}_{\Gamma}=\bm{u}_{\Gamma}^i.
\end{eqnarray*}


Define the local Schur complement matrix $S_i$ and the corresponding right-hand side vector $F_i$ as
\begin{eqnarray*}
	S_i=A_{\Gamma\Gamma}^i-A_{\Gamma I}^i(A_{II}^i)^{-1}A_{I\Gamma}^i,\
	F_i=\bm{f}_{\Gamma}^i-A_{\Gamma I}^i(A_{II}^i)^{-1}\bm{f}_{I}^i.
\end{eqnarray*}
Then the global interface unknowns $\bm{u}_{\Gamma}$ can be obtained by solving the global interface problem
\begin{eqnarray}
	\label{interface-d}
	S\bm{u}_{\Gamma}=F,
\end{eqnarray}
where
\begin{eqnarray*}
	S:=\sum_{i=1}^{N_s}R_i^{T}S_iR_i,\ F:=\sum_{i=1}^{N_s}R_i^{T}F_i,
\end{eqnarray*}
and the interior unknowns in each subdomain $D_i$ can be obtained through (\ref{eq-interior}).

It is worth noting that the interface system described above exhibits a smaller but
denser structure compared to the original global system. Furthermore, the condition number of the Schur complement matrix $S$ is generally better than the original global stiffness matrix \cite{T1994ddm}.  Notably, the computational cost associated
with calculating $S$ and $F$  primarily arises from the inversion of $A_{II}^i$, for $ 1\leq i \leq N_s$, particularly when dealing with a large number of degrees of freedom.  In this study,
our objective is to explore the application of the domain decomposition method to
stochastic partial differential equations (PDEs). In this context, all the { aforementioned} matrices are influenced by random inputs, introducing significant challenges
when attempting to solve the stochastic interface system.

\section{Stochastic domain decomposition based on Variable-separation method}
\label{sec-sddm}
In this section, we will present an offline-online method for the interface problem of the SPDEs, which builds a relation between the random inputs and the stochastic interface solution. The purpose of the offline stage is to construct all the components that are needed in the online stage. Details about the offline stage are provided in Subsections \ref{sec-sf}-\ref{sec-reduce model}. Assuming that the original domain $D$ is divided into $N_s$ non-overlapping subdomains, the interface problem of stochastic PDEs can
be expressed as follows  (referring to equation (\ref{interface-d})):
\begin{eqnarray}
	\label{eq-interface}
	S(\bm{\xi})\bm{u}_{\Gamma}(\bm{\xi})=F(\bm{\xi}),
\end{eqnarray}
where
\begin{eqnarray*}
	\label{sub-global}
	S(\bm{\xi}):=\sum_{i=1}^{N_s}R_i^{T}S_i(\bm{\xi})R_i,\ F(\bm{\xi}):=\sum_{i=1}^{N_s}R_i^{T}F_i(\bm{\xi}),
\end{eqnarray*}
and
\begin{align}
	\label{eq-S-F-1}
	S_i(\bm{\xi})&=A_{\Gamma\Gamma}^i(\bm{\xi})-A_{\Gamma I}^i(\bm{\xi})(A_{II}^i(\bm{\xi}))^{-1}A_{I\Gamma}^i(\bm{\xi}),\\
	\label{eq-S-F-2}
	F_i(\bm{\xi})&=\bm{f}_{\Gamma}^i(\bm{\xi})-A_{\Gamma I}^i(\bm{\xi})(A_{II}^i(\bm{\xi}))^{-1}\bm{f}_{I}^i(\bm{\xi}).
\end{align}

For an elementary exposition of the methodology, we consider dividing $D$ into two subdomains $D_1$ and $D_2$, and redefine the interface of $D_1$ and $D_2$ by $\Gamma$ instead of $\Gamma_{12}$, as shown in Figure \ref{ddm}. It should be clarified that all the discussions in two subdomains can be generalized to the case of multiple subdomains. Then we have the interface problem of stochastic PDEs as follows
\begin{eqnarray}
	\label{eq-spde-interface}
	(S_1(\bm{\xi})+S_2(\bm{\xi}))\bm{u}_{\Gamma}(\bm{\xi})=F_1(\bm{\xi})+F_2(\bm{\xi}).
\end{eqnarray}
Note that, for each sample $\bm{\xi}\in \Omega$, the calculation of $S_1(\bm{\xi}),S_2(\bm{\xi}),F_1(\bm{\xi}),F_2(\bm{\xi})$ by equations (\ref{eq-S-F-1}-\ref{eq-S-F-2}) depends on the inversion of the full order matrices $A_{II}^1(\bm{\xi}), A_{II}^2(\bm{\xi})$, which will substantially impact on the computation efficiency of solving the interface problem (\ref{eq-spde-interface}).
To improve computational efficiency, we aim to reconstruct equation (\ref{eq-spde-interface}) such that
$S(\bm{\xi})$ and $F(\bm{\xi})$ are affine with respect to $\bm{\xi}$, i.e.
\begin{eqnarray}
	\label{eq-S-F-Aff1}
	S(\bm{\xi})
	=\sum_{j=1}^{m_S}\hat{\eta}^j(\bm{\xi})\mathcal{\hat{X}}^j,\\
	\label{eq-S-F-Aff2}
	F(\bm{\xi})
	=\sum_{j=1}^{m_F}\hat{\gamma}^j(\bm{\xi}){\hat{F}}^j.
\end{eqnarray}
The details of the process and definitions for $\hat{\eta}^j(\bm{\xi}), \hat{\gamma}^j(\bm{\xi})$, $\mathcal{\hat{X}}^j$ and ${\hat{F}}^j$ will be introduced in Subsection \ref{sec-sf}. By equations (\ref{eq-S-F-Aff1}-\ref{eq-S-F-Aff2}), we have the stochastic interface problem such as
\begin{eqnarray}
	\label{eq-Intaff}
	\sum_{j=1}^{m_S}\hat{\eta}^j(\bm{\xi})\mathcal{\hat{X}}^j \bm{u}_{\Gamma}(\bm{\xi})=\sum_{j=1}^{m_F}\hat{\gamma}^j(\bm{\xi}){\hat{F}}^j.
\end{eqnarray}
The matrices $\mathcal{\hat{X}}^j$ and ${\hat{F}}^j$ are independent of the random variables $\bm{\xi}$, and their computation is a one-time operation.
For any $\bm{\xi}\in \Omega$, we just need to solve equation (\ref{eq-Intaff}) instead of equation (\ref{eq-interface}) involving the inversion of the full order matrices.

\subsection{Assemble strategies for S and F}
\label{sec-sf}
Now we describe the detail of the strategy for constructing equations (\ref{eq-S-F-Aff1}-\ref{eq-S-F-Aff2}).
With the assumption (\ref{eq-vs-form}) of affine decomposition,
the stochastic matrices in equations (\ref{eq-S-F-1}-\ref{eq-S-F-2}) can be written as follows
\begin{equation}
	\label{eq-affmatrix}
	\begin{aligned}
		A_{II}^{i}(\bm{\xi})&=\sum_{j=1}^{m_{a_i}}p^{ij}(\bm{\xi})A_{II}^{ij},
		~A_{I\Gamma}^{i}(\bm{\xi})=\sum_{j=1}^{m_{a_i}}p^{ij}(\bm{\xi})A_{I\Gamma}^{ij},\\
		A_{\Gamma I}^{i}(\bm{\xi})&=\sum_{j=1}^{m_{a_i}}p^{ij}(\bm{\xi})A_{\Gamma I}^{ij},
		~A_{\Gamma\Gamma}^{i}(\bm{\xi})=\sum_{j=1}^{m_{a_i}}p^{ij}(\bm{\xi})A_{\Gamma\Gamma}^{ij},\\
	\end{aligned}
\end{equation}
where $A_{II}^{ij},A_{I\Gamma}^{ij},A_{\Gamma I}^{ij},A_{\Gamma\Gamma}^{ij}$ are independent of random variables $\bm{\xi}$. Thus the first item of $S_i(\bm{\xi})$ in equation (\ref{eq-S-F-1}), i.e., $A_{\Gamma\Gamma}^{i}(\bm{\xi})$ is affine with respect to $\bm{\xi}$ naturally. To achieve affine expression for the second item of $S_i(\bm{\xi})$, i.e., $A_{\Gamma I}^i(\bm{\xi})(A_{II}^i(\bm{\xi}))^{-1}A_{I\Gamma}^i(\bm{\xi})$, we perform it as follows.

$\bullet$ \textit{Step 1: Construct the low-rank approximation of $X(\bm{\xi})=(A_{II}^i(\bm{\xi}))^{-1}A_{I\Gamma}^i(\bm{\xi})$ such as
	\begin{eqnarray} 	
		\label{eq-X-appr}
		X_{N}(\bm{\xi}):=\sum_{j=1}^{N_{S_i}}\beta_j(\bm{\xi})X_j.
\end{eqnarray}}

First, we rewrite $(A_{II}^i(\bm{\xi}))^{-1}A_{I\Gamma}^i(\bm{\xi})$ as the following stochastic algebraic system
\begin{eqnarray} 	
	\label{eq-stepas}
	A_{II}^{i}(\bm{\xi})X(\bm{\xi})=A_{I\Gamma}^{i}(\bm{\xi}),
\end{eqnarray}
with assumptions of affine decomposition (\ref{eq-affmatrix}), we have
\begin{eqnarray*}
	\sum_{j=1}^{m_{a_i}}p^{ij}(\bm{\xi})A_{II}^{ij}X(\bm{\xi})=\sum_{j=1}^{m_{a_i}}p^{ij}(\bm{\xi})A_{I\Gamma}^{ij}.
\end{eqnarray*}
Let $n_{\Gamma}$ be the number of the interface unknowns, we rewrite the stochastic matrix as $X(\bm{\xi})=[x_1^{i}(\bm{\xi}), x_2^{i}(\bm{\xi}), \cdots,x_{n_{\Gamma}}^{i}(\bm{\xi})]$. It follows that to get the solution of equation (\ref{eq-stepas}) is equivalent to solving the following $n_{\Gamma}$ equations
\begin{eqnarray} 	
	\label{eq-stepass}
	\sum_{j=1}^{m_{a_i}}p^{ij}(\bm{\xi})A_{II}^{ij} x_k^{i}(\bm{\xi})=\sum_{j=1}^{m_{a_i}}p^{ij}(\bm{\xi})\alpha_k^{ij}, ~~ 1\leq k\leq n_{\Gamma},
\end{eqnarray}
where $\alpha_k^{ij}$ represents the $k$-th column of the matrix $A_{I\Gamma}^{ij}$ for  $1\leq k\leq n_{\Gamma}$, and $1\leq j\leq m_a$.
For each equation of (\ref{eq-stepass}), we adopt Algorithm \ref{algorithm-vs-sSAS} proposed in Subsection \ref{sec-vs} to get the low-rank approximation of $x_k^{i}(\bm{\xi})$ such as (\ref{eq-approx}). Consequently, the low-rank approximation of $X(\bm{\xi})$ in the form of (\ref{eq-X-appr}) can be given by rearranging the low-rank approximations of $x_k^{i}(\bm{\xi})$ for $1\leq k\leq n_{\Gamma}$.

$\bullet$ \textit{Step 2: Assemble affine expression for
	\begin{eqnarray}
		\label{xdef} 		
		\mathcal{X}^i(\bm{\xi})=-A_{\Gamma I}^i(\bm{\xi})(A_{II}^i(\bm{\xi}))^{-1}A_{I\Gamma}^i(\bm{\xi}).
	\end{eqnarray}
}

Based on the low-rank representation of $X(\bm{\xi})$ and equation (\ref{eq-X-appr}), we have
\begin{eqnarray} 	
	\label{eq-stepbs}
	\mathcal{X}^i(\bm{\xi})\approx -A_{\Gamma I}^i(\bm{\xi})\sum_{k=1}^{N_{S_i}}\beta_k(\bm{\xi})X_k
	=-\sum_{j=1}^{m_{a_i}}\sum_{k=1}^{N_{S_i}}p^{ij}(\bm{\xi})\beta_k(\bm{\xi})A_{\Gamma I}^{ij}X_k.
\end{eqnarray}
The affine decomposition of $A_{\Gamma I}^i(\bm{\xi})$ in equation (\ref{eq-affmatrix}) yields the second equality.

To simplify notation, the affine expression for $\mathcal{X}^i(\bm{\xi})$  can be obtained by using the single-index notation as follows
\begin{eqnarray} 	
	\label{eq-stepb}
	\mathcal{X}^i(\bm{\xi})=\sum_{j\in J}\eta_{ij}(\bm{\xi})\mathcal{X}^{ij},
\end{eqnarray}
where $J= \{1,2,\cdots,m_{a_i}\}\times\{1,2,\cdots,N_{S_i}\}$,  $\eta_{ij}(\bm{\xi})=p^{ij_1}(\bm{\xi})\beta_{j_2}(\bm{\xi})$, and $\mathcal{X}^{ij}=-A_{\Gamma I}^{ij_1}X_{j_2}$ is matrix of $n_{\Gamma} \times n_{\Gamma}$ for arbitrary $1\leq j\leq m_{a_i}N_{S_i}$.
Subsequently, we have the low-rank representation for $S_i(\bm{\xi})$ based on equation (\ref{eq-S-F-1}) as follows
\begin{eqnarray*}
	S_i(\bm{\xi})= A_{\Gamma\Gamma}^i(\bm{\xi})+\mathcal{X}^i(\bm{\xi})
	\approx\sum_{j=1}^{m_{a_i}}p^{ij}(\bm{\xi})A_{\Gamma\Gamma}^{ij}+\sum_{j=1}^{{m_{a_i}}N_{S_i}}\eta_{ij}(\bm{\xi})\mathcal{X}^{ij},
\end{eqnarray*}
where the first equality follows from the definition of $\mathcal{X}^{i}(\bm{\xi})$ in equation (\ref{xdef}),  the second equality follows from the affine decomposition of $A_{\Gamma \Gamma}^i(\bm{\xi})$ and equation (\ref{eq-stepb}). Consequently, we have the low-rank representation for $S(\bm{\xi})$ as follows
\begin{equation}
	\label{eq-S-repre}
	\begin{aligned}
		S(\bm{\xi})&= S_1(\bm{\xi})+S_2(\bm{\xi})\\
		&\approx\sum_{j=1}^{{m_{a_1}}}p^{1j}(\bm{\xi})A_{\Gamma\Gamma}^{1j}+\sum_{j=1}^{{m_{a_2}}}p^{2j}(\bm{\xi})A_{\Gamma\Gamma}^{2j}+\sum_{j=1}^{m_{a_1}N_{S_1}}\eta_{1j}(\bm{\xi})\mathcal{X}^{1j}+\sum_{j=1}^{m_{a_2}N_{S_2}}\eta_{2j}(\bm{\xi})\mathcal{X}^{2j}\\
		&=\sum_{j=1}^{m_{a_1}(N_{S_1}+1)+m_{a_2}(N_{S_2}+1)}\hat{\eta}^j(\bm{\xi})\mathcal{\hat{X}}^j,
	\end{aligned}
\end{equation}
where the last equality comes from stacking the variables and sorting the corresponding indices via a single one.
Algorithm \ref{algorithm-S} outlines the assembling process of $S(\bm{\xi})$.

\begin{algorithm}
	\caption{The Assemble Process of $S(\bm{\xi})$}
	\textbf{Input}: The stochastic matrices $A_{II}^{i}(\bm{\xi}),A_{I\Gamma}^{i}(\bm{\xi}),A_{\Gamma I}^{i}(\bm{\xi}),A_{\Gamma\Gamma}^{i}(\bm{\xi})$, $i=1,2$.  \\
	\textbf{Output}: The low-rank representation $S(\bm{\xi})=\sum_{j=1}^{m_{a_1}(N_{S_1}+1)+m_{a_2}(N_{S_2}+1)}\hat{\eta}^j(\bm{\xi})\mathcal{\hat{X}}^j$.\\
	~1:~~Get the approximation of $X(\bm{\xi})=(A_{II}^i(\bm{\xi}))^{-1}A_{I\Gamma}^i(\bm{\xi})$ by solving equations (\ref{eq-stepass})
	
	~~~~~use Algorithm \ref{algorithm-vs-sSAS};\\
	~2:~~Assemble the affine expression of $-A_{\Gamma I}^i(\bm{\xi})(A_{II}^i(\bm{\xi}))^{-1}A_{I\Gamma}^i(\bm{\xi})$ by (\ref{eq-stepbs});\\
	~3:~~Assemble $S_i(\bm{\xi})$ based on (\ref{eq-affmatrix}) and the expression derived in step 2;\\
	~4:~~Assemble the low-rank representation of $S(\bm{\xi})$ by (\ref{eq-S-repre}).
	\label{algorithm-S}
\end{algorithm}

Similarly, with the assumption (\ref{eq-vs-form}) of affine decomposition
\begin{eqnarray}
	\label{eq-affvector}
	\bm{f}_I^{i}(\bm{\xi})=\sum_{j=1}^{m_{b_i}}q^{ij}(\bm{\xi})\bm{f}_I^{ij},~~ \bm{f}_{\Gamma}^{i}(\bm{\xi})=\sum_{j=1}^{m_{b_i}}q^{ij}(\bm{\xi})\bm{f}_{\Gamma}^{ij},
\end{eqnarray}
where $\bm{f}_I^{ij}, \bm{f}_{\Gamma}^{ij}$ are independent of $\bm{\xi}$. Then following the assemble strategy of $S(\bm{\xi})$, the low-rank representation for $F(\bm{\xi})$ can be constructed as follows
\begin{eqnarray}
	\label{eq-F-repre}
	\begin{aligned}
		F(\bm{\xi})&= F_1(\bm{\xi})+F_2(\bm{\xi})\\
		&\approx\sum_{j=1}^{m_{b_1}}q^{1j}(\bm{\xi})\bm{f}_{\Gamma}^{1j}+\sum_{j=1}^{m_{b_2}}q^{2j}(\bm{\xi})\bm{f}_{\Gamma}^{2j}+\sum_{j=1}^{m_{a_1}N_{F_1}}\gamma_{1j}(\bm{\xi}){F}^{1j}+\sum_{j=1}^{m_{a_2}N_{F_2}}\gamma_{2j}(\bm{\xi}){F}^{2j}\\
		&=\sum_{j=1}^{m_{b_1}+m_{b_2}+m_{a_1}N_{F_1}+m_{a_2}N_{F_2}}\hat{\gamma}^j(\bm{\xi}){\hat{F}}^j.
	\end{aligned}
\end{eqnarray}

We present the detail of the process of assembling $F(\bm{\xi})$ in Algorithm \ref{algorithm-F}.

\begin{algorithm}
	\caption{The Assemble Process of $F(\bm{\xi})$ }
	\textbf{Input}: The stochastic matrices $A_{II}^{i}(\bm{\xi}),A_{\Gamma I}^{i}(\bm{\xi})$ and vectors $\bm{f}_I^{i}(\bm{\xi}),\bm{f}_{\Gamma}^{i}(\bm{\xi})$, $i=1,2$.  \\
	\textbf{Output}: The low-rank representation $	F(\bm{\xi})=\sum_{j=1}^{m_{b_1}+m_{b_2}+m_{a_1}N_{F_1}+m_{a_2}N_{F_2}}\hat{\gamma}^j(\bm{\xi}){\hat{F}}^j$.\\
	~1:~~Get the approximation of $X(\bm{\xi})=(A_{II}^i(\bm{\xi}))^{-1}\bm{f}_{I}^{i}(\bm{\xi})$ by Algorithm \ref{algorithm-vs-sSAS};\\
	~2:~~Assemble the affine expression of $-A_{\Gamma I}^i(\bm{\xi})(A_{II}^i(\bm{\xi}))^{-1}\bm{f}_{I}^{i}(\bm{\xi})$;\\
	~3:~~Assemble $F_i(\bm{\xi})$ based on (\ref{eq-affvector}) and the expression derived in step 2;\\
	~4:~~Assemble the low-rank representation of $F(\bm{\xi})$ by (\ref{eq-F-repre}).
	\label{algorithm-F}
\end{algorithm}

\begin{remark}
	Note that $m_{a_1}(N_{S_1}+1)+m_{a_2}(N_{S_2}+1)$ is the maximum number of the separated terms of $S(\bm{\xi})$. For practical problem, the number of the separated terms of $S(\bm{\xi})$ is smaller than $m_{a_1}(N_{S_1}+1)+m_{a_2}(N_{S_2}+1)$ when there exist $p^{1j}(\bm{\xi})=p^{2k}(\bm{\xi})$, $j=1,2,\cdots,m_{a_1}, k=1,2,\cdots,m_{a_2}$, the same is true for $F(\bm{\xi})$.
\end{remark}

\subsection{Reduced model representation for the stochastic interface problem}
\label{sec-reduce model}
As we mentioned before, once the separation approximations of $S(\bm{\xi})$ and $F(\bm{\xi})$ are available, the stochastic interface problem (\ref{eq-spde-interface}) can be given by
\begin{eqnarray}
	\label{eq-interf}
	\sum_{j=1}^{m_S}\hat{\eta}^j(\bm{\xi})\mathcal{\hat{X}}^j \bm{u}_{\Gamma}(\bm{\xi})=\sum_{j=1}^{m_F}\hat{\gamma}^j(\bm{\xi}){\hat{F}}^j, \ \forall ~\bm{\xi} \in \Omega,
\end{eqnarray}
where $m_S=m_{a_1}(N_{S_1}+1)+m_{a_2}(N_{S_2}+1)$ and $m_F=m_{b_1}+m_{b_2}+m_{a_1}N_{F_1}+m_{a_2}N_{F_2}$. The amended model defined in (\ref{eq-interf}) for the interface problem is a linear algebraic system with $n_\Gamma$ unknowns. Although the amended model needs much less computation effort compared with the original model (\ref{eq-spde-interface}), it may be not a very small-scale problem because the amended model defined in (\ref{eq-interf}) involves the discrete degree of the original full model. {In order to significantly improve the computation efficiency, we want to get the reduced model representation of the stochastic interface problem, i.e.,
	\begin{eqnarray}
		\label{eq-interu}
		\bm{u}_{\Gamma}(\bm{\xi})\approx\sum_{i=1}^{M}\zeta_i(\bm{\xi})\bm{c}_{i},
	\end{eqnarray}
	where $M$ is the number of the separated terms of $\bm{u}_{\Gamma}(\bm{\xi})$, each $\bm{c}_i,i=1,\cdots,M$ is a vector of $n_\Gamma$-dimension. Then the functional decomposition expression
	of the stochastic solution for the interface problem can be written as follows
	\begin{eqnarray}
		\label{eq-interu1}
		{u}_{\Gamma}(\bm{\xi})\approx\sum_{i=1}^{M}\zeta_i(\bm{\xi}){c}_{i}(x),
	\end{eqnarray}
	with ${c}_{i}(x)=\sum_{j=1}^{n_\Gamma}(\bm{c}_{i})_j\psi_{j}(x),i=1,\cdots,M$ and $\{\psi_{j}\}_{j=1}^{n_\Gamma}$ being the corresponding basis functions of the FE space $\mathcal{V}_h$ on the interface $\Gamma$.}
We call this reduced model representation, which is more applicable to solving the subproblems on subdomains.

We employ the VS method for stochastic algebraic systems presented in Subsection \ref{sec-vs} to derive the reduced model representation (\ref{eq-interu}), which is beneficial to construct the efficient surrogate model of the subproblems.

The online stage of the SDD-VS method is to use the output of the offline stage to recover the solution to the stochastic interface problem for a large number of new samples. The online stage is efficient thanks to the reduced model representation (\ref{eq-interu}) and (\ref{eq-interu1}) for the stochastic interface solution.

{ \subsection{The numerical method for the subproblems}}
\label{sec-subproblems}
In the above subsections, we have a detailed description of the SDD-VS method for the stochastic interface problem. Once we get the stochastic solution ${u}_{\Gamma}^i(\bm{\xi})$ of the interface, the interior solution ${u}_{I}^i(\bm{\xi})$ in the subdomain $D_i$ can be obtained by solving the following stochastic problem
\begin{eqnarray}
	\label{eq-Sub-SPDE}
	\left\{
	\begin{aligned}
		\mathcal{L}(x,{\xi};{u}_{I}^i(\bm{\xi}))&=f(x,\bm{\xi}), ~\ \forall ~x \in D_i, ~ \bm{\xi} \in \Omega,\\	
		\mathcal{B}(x,{\xi};{u}_{I}^i(\bm{\xi}))&=
		g(x,\bm{\xi}), ~\ \forall ~x \in \partial D_i\backslash(\partial D_i\cap \Gamma), ~ \bm{\xi} \in \Omega,\\
		\mathcal{B}(x,{\xi};{u}_{I}^i(\bm{\xi}))&=h_i(x,\bm{\xi}), \ \forall ~x \in \partial D_i\cap \Gamma, ~ \bm{\xi} \in \Omega,
	\end{aligned}
	\right.
\end{eqnarray}
where $h_i(x,\bm{\xi})$ is defined by the interface solution ${u}_{\Gamma}^i(\bm{\xi})$ such as equation (\ref{eq-interu1}) on the interface boundary, i.e., $\partial D_i\cap \Gamma$.
To address the computational complexity associated with solving the stochastic partial differential equation (\ref{eq-Sub-SPDE}), we attempt to apply model reduction methods to construct an efficient surrogate model.  Model reduction methods have been proposed to reduce the computation complexity especially when the full model are expensive to perform numerical simulations. These methods construct an approximate model with lower dimensionality but still describe important aspects of the full model. The reduced basis method is one of the model order reduction methods and usually provides an efficient and reliable approximation of the input-output relationship \cite{Canuto2009a, Elman2013reduced, Chen2013a, Drohmann2012Reduced, Jiang2017Model, Quarteroni2015Reduced}. Another class of model reduction methods is based on the variable-separation (VS) method. As an example, the proper generalized decomposition (PGD) method has been used in solving stochastic partial differential equations (SPDEs) \cite{ cho2016numerical}. The PGD method constructs optimal reduced basis from a double orthogonality criterium  \cite{nouy2009generalized,nouy2010proper}, and it requires the solutions of a few uncoupled deterministic problems solved by classical deterministic solution techniques and the solutions of stochastic algebraic are solved by classical spectral stochastic methods. We note that PGD requires many iterations with the arbitrary initial guess to compute each term in the separated expansion at each enrichment step. This will deteriorate the simulation efficiency.

Here, we would like to adopt the VS method we proposed in \cite{Li2017a} for the stochastic problem  (\ref{eq-Sub-SPDE}) to get a separated representation for the solution without iterations at each enrichment step. Moreover, the proposed VS method can alleviate the ``curse of dimensionality'' when dealing with problems in high-dimensional stochastic spaces.

Finally, we summarize the SDD-VS method for stochastic PDEs in Figure \ref{flowchart-SDD-VS}.
\tikzset{meta box/.style={draw,black},
	punkt/.style={meta box,rectangle,rounded corners,minimum height=2em,minimum width=6em,text width=25em},
	punktd/.style={meta box,rectangle,rounded corners,minimum height=2em,minimum width=6em,text width=17em}}
\tikzset{add dimen/.code 2 args={\pgfkeysgetvalue{/pgf/minimum #1}\tikz@dimen@min
		\expandafter\tikz@expand@dimen\expandafter{\tikz@dimen@min + #2 * 2em}{#1}},
	wider/.style={add dimen={width}{#1}},	higher/.style={add dimen={height}{#1}},}

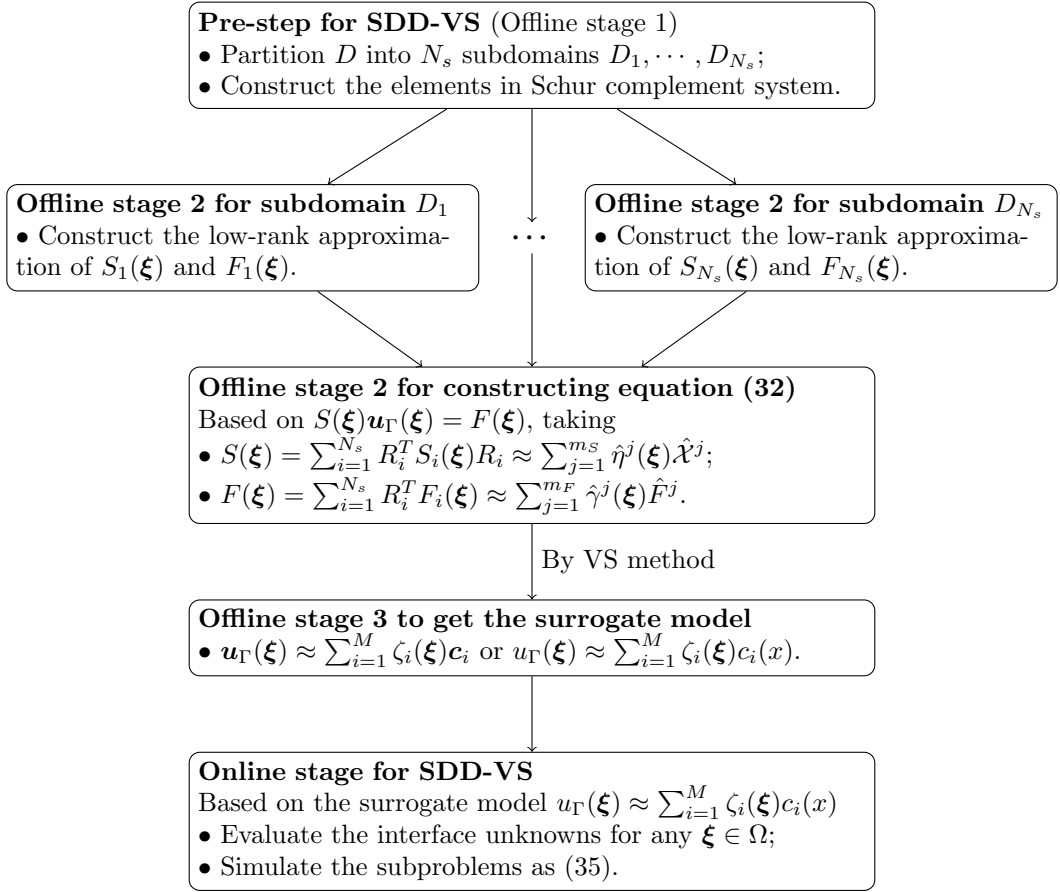
\begin{figure}[htbp]
	\centering
	\begin{tikzpicture}[scale=1]
		\node[punkt, wider=3,higher=2.5] (pre)
		{\textbf{Pre-step for SDD-VS} (Offline stage 1)\\
			$\bullet$ Partition $D$ into $N_s$ subdomains $D_1, \cdots, D_{N_s}$;\\
			$\bullet$ Construct the elements in Schur complement system.};
		
		\node[draw=white, below=1.5 of pre] (Di)
		{$\bm{\cdots}$};
		
		\node[punktd, wider=3, higher=2.5, left=0.3 of Di] (D1)
		{\textbf{Offline stage 2 for subdomain $D_1$}\\ $\bullet$ Construct the low-rank approximation of $S_{1}(\bm{\xi})$ and $F_{1}(\bm{\xi})$.};
		
		\node[punktd, wider=3, higher=2.5, right=0.3 of Di] (Dn)
		{\textbf{Offline stage 2 for subdomain $D_{N_s}$}\\ $\bullet$ Construct the low-rank approximation of $S_{N_s}(\bm{\xi})$ and $F_{N_s}(\bm{\xi})$.};
		
		\node[punkt, wider=3, higher=2.5, below=1.5 of Di] (inter)
		{\textbf{Offline stage 2 for constructing equation (\ref{eq-interf})} \\
			Based on $S(\bm{\xi})\bm{u}_{\Gamma}(\bm{\xi})=F(\bm{\xi})$, taking\\
			$\bullet$ $S(\bm{\xi})=\sum_{i=1}^{N_s}R_i^{T}S_i(\bm{\xi})R_i\approx\sum_{j=1}^{m_S}\hat{\eta}^j(\bm{\xi})\mathcal{\hat{X}}^j $;\\
			$\bullet$  $F(\bm{\xi})=\sum_{i=1}^{N_s}R_i^{T}F_i(\bm{\xi})\approx\sum_{j=1}^{m_F}\hat{\gamma}^j(\bm{\xi}){\hat{F}}^j$.};
		
		\node[punkt, wider=3,higher=2.5, below=1 of inter] (solu)
		{\textbf{Offline stage 3 to get the surrogate model}\\
			$\bullet$ 	$\bm{u}_{\Gamma}(\bm{\xi})\approx\sum_{i=1}^{M}\zeta_i(\bm{\xi})\bm{c}_{i}$ or ${u}_{\Gamma}(\bm{\xi})\approx\sum_{i=1}^{M}\zeta_i(\bm{\xi})c_{i}(x)$.};

		\node[punkt, wider=3,higher=2.5, below=1 of solu] (online)
		{\textbf{Online stage for SDD-VS}\\
			Based on the surrogate model ${u}_{\Gamma}(\bm{\xi})\approx\sum_{i=1}^{M}\zeta_i(\bm{\xi})c_{i}(x)$\\
			$\bullet$  Evaluate the interface unknowns for any $\bm{\xi}\in\Omega$; \\
			$\bullet$  {Simulate the subproblems as (\ref{eq-Sub-SPDE})}.};	
		\graph{
			(pre)->(D1)->(inter);
			(pre)->(Di)->(inter)->["By VS method"right](solu)-> (online);
			(pre)->(Dn)->(inter);
		};		
	\end{tikzpicture}
	\caption{Flowchart of the SDD-VS method.}
	\label{flowchart-SDD-VS}
\end{figure}

\section{Numerical experiments}
\label{sec-numerical examples}
In this section, we will present various numerical results to demonstrate the applicability and efficiency of the proposed SDD-VS method on several stochastic PDEs. For each problem, we seek a separate representation approximation to the interface problem induced by the model depending on random variables.
All the numerical experiments in this paper were run in Python on a Dell desktop with Intel Core i7-4970 CPU @3.60GHz and 16GB of RAM.
In Section \ref{num1}, we consider the one-dimensional (1D) stochastic diffusion equation to illustrate the performance of the proposed SDD-VS method. In Section \ref{num2}, we study the SDD-VS method for a two-dimensional (2D) stochastic diffusion equation with three subdomains. The 2D stochastic convection-diffusion equation with high-dimensional random inputs is considered in Section \ref{num3}.

In order to quantify the accuracy of the proposed SDD-VS method, we use the relative mean error $\epsilon$ for the stochastic interface problem as follows
\begin{eqnarray}
	\label{eq-meanerror}
	\epsilon=\frac{1}{N}\sum_{i=1}^{N}\frac{|u_{\Gamma}(\bm{\xi}_i)-\hat{u}_{\Gamma}(\bm{\xi}_i)|}{|u_{\Gamma}(\bm{\xi}_i)|},
\end{eqnarray}
where $N$ is the number of samples used to compute the mean error, $\hat{u}_{\Gamma}(\bm{\xi})$ is the approximation of the interface unknowns obtained by the SDD-VS method, $u_{\Gamma}(\bm{\xi})$ is the restriction of the reference solution, which is solved by the finite element method over the whole domain $D$, on the interface part $\Gamma$.

\subsection{1D stochastic diffusion equation}
\label{num1}
We begin by considering a 1D stochastic diffusion equation with homogeneous Dirichlet boundary conditions
\begin{eqnarray*}
	\left\{
	\begin{aligned}
		-\frac{d}{dx}(c(x;\bm{\xi}) \frac{du(x;\bm{\xi})}{dx})&=f(x;\bm{\xi}),  \ \forall ~x\in D, \\	
		u(0;\bm{\xi})=u(1;\bm{\xi})&=0,	
	\end{aligned}
	\right.
\end{eqnarray*}
where the original domain $D=[0,1]$ is divided into two subdomains $D_1=[0,0.5]$, $D_2=\left(0.5,1\right]$, and the random coefficient $c(x;\bm{\xi})$, the source function $f(x;\bm{\xi})$ are defined as follows
\begin{eqnarray*}
	\begin{aligned}
		c(x;\bm{\xi})&=
		\left\{
		\begin{aligned}
			&\bm{\xi}x+4, ~~&\text{when} ~~x\in D_1,\\
			&x+1, &\text{when} ~~x\in D_2,\\
		\end{aligned}
		\right.
		\\
		f(x;\bm{\xi})&=
		\left\{
		\begin{aligned}
			&\cos(2\pi x), &\text{when} ~~x\in D_1,\\
			&\bm{\xi}^2x, &\text{when} ~~x\in D_2.\\
		\end{aligned}
		\right.
	\end{aligned}	
\end{eqnarray*}

In this example, we set the random variable $\bm{\xi}$ to be a truncated Gaussian distribution with mean 0, standard deviation 1, and range $[-3,3]$. Here the reference solution is calculated by the finite element method with mesh size $h=1/1000$. For the stochastic problems of one dimension, the reduced model (\ref{eq-interf})
\begin{eqnarray*}
	\sum_{j=1}^{m_S}\hat{\eta}^j(\bm{\xi})\mathcal{\hat{X}}^j \bm{u}_{\Gamma}(\bm{\xi})=\sum_{j=1}^{m_F}\hat{\gamma}^j(\bm{\xi}){\hat{F}}^j, \ \forall ~\bm{\xi} \in \Omega,
\end{eqnarray*}
is just a linear algebraic equation, where  $\{\mathcal{\hat{X}}^j\}_{j=1}^{m_S}$ and $\{\hat{F}^j\}_{j=1}^{m_F}$ are constants. Thus we can obtain the simple analytic expression for the stochastic interface unknowns
\begin{eqnarray}
	\label{eq-intso1}
	\bm{u}_{\Gamma}(\bm{\xi})=\dfrac{\sum_{j=1}^{m_F}\hat{\gamma}^j(\bm{\xi}){\hat{F}}^j}{	\sum_{j=1}^{m_S}\hat{\eta}^j(\bm{\xi})\mathcal{\hat{X}}^j}, \ \forall ~\bm{\xi} \in \Omega,
\end{eqnarray}
without further calculations. Therefore, most of the cost in the offline phase concentrates on the assemble process of $S_i(\bm{\xi}), F_i(\bm{\xi})$, $i=1,2$ in equations (\ref{eq-S-F-1}-\ref{eq-S-F-2}), i.e. the construction of their low-rank representation.
To this end, we choose $|\Xi|=20$ samples for training in Algorithm \ref{algorithm-vs-sSAS}.

\begin{figure}[htbp]
	\centering
	\includegraphics[width=3.8in, height=2.7in]{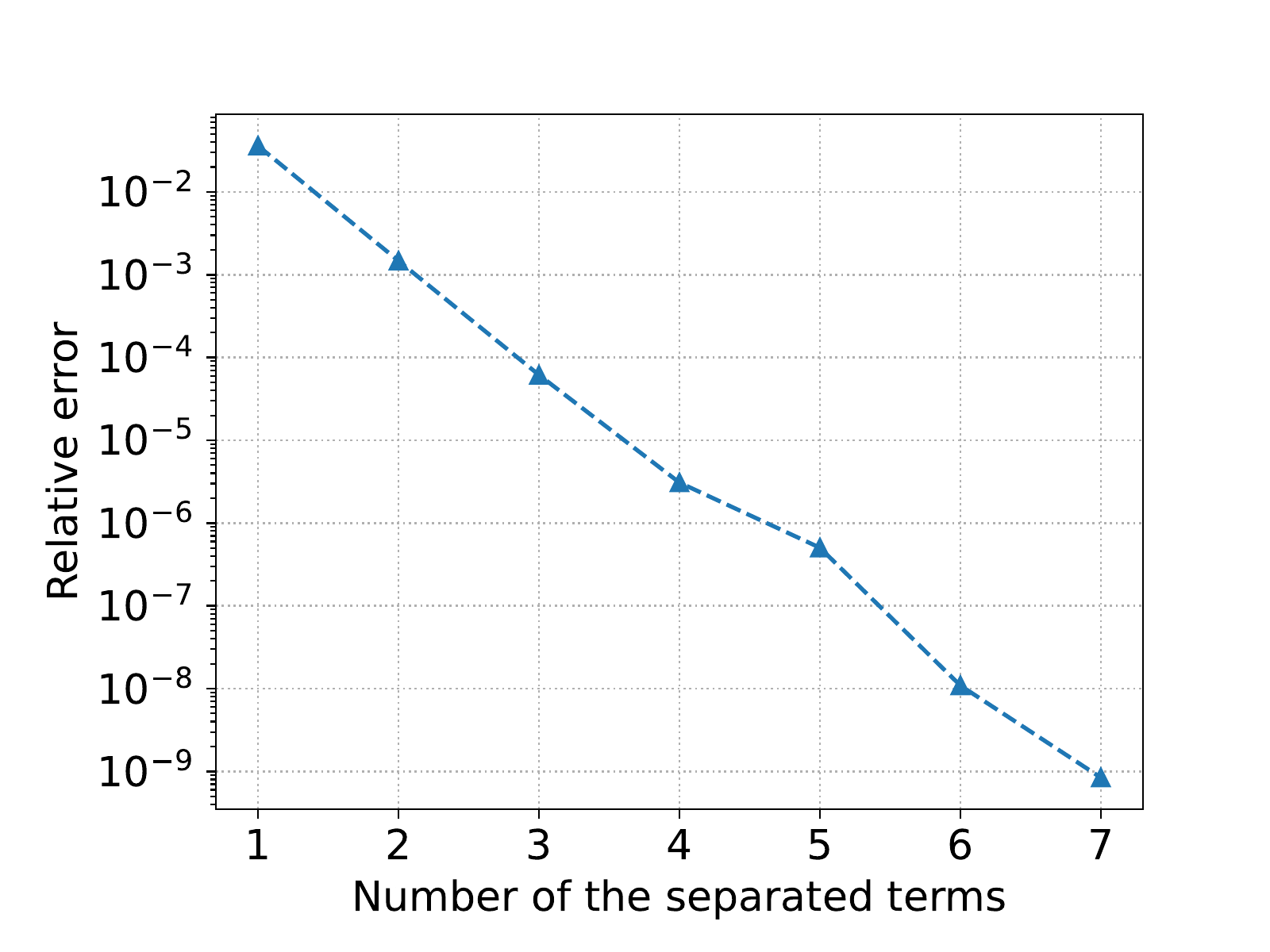}
	\caption{Comparison of the relative mean error corresponding to the different numbers of the separated terms $N_{S_1}$.}
	\label{fig2}
\end{figure}

It should be noted that the relative mean error in all the assemble processes becomes smaller as the number of the separated terms increases.  We take the construction of $S_1(\bm{\xi})$ as an example, and depict the relative mean error versus the number of the separated terms for the VS method in Figure \ref{fig2}, where the relative mean error $\epsilon$ is calculated by equation (\ref{eq-meanerror}) with $N=10^4$ test samples.
This could provide evidence for the selection of $N_{S_1}=4$ such that the relative mean error $\epsilon < 10^{-5}$ is small enough. Similarly, we take  $ N_{S_2}=1,  N_{F_1}=4, N_{F_2}=1 $ for the construction of $S_2(\bm{\xi})$, $F_1(\bm{\xi})$ and $F_2(\bm{\xi})$ respectively. Thus we have  $m_S=m_{a_1}(N_{S_1}+1)+m_{a_2}(N_{S_2}+1)=12$ and $m_F=m_{b_1}+m_{b_2}+m_{a_1}N_{F_1}+m_{a_2}N_{F_2}=11$ in the reduced model (\ref{eq-interf}), where $m_{a_1}=2, m_{a_2}=1$ and $m_{b_1}=m_{b_2}=1$.

\begin{figure}[htb]
	\centering
	\includegraphics[width=3.8in, height=2.7in]{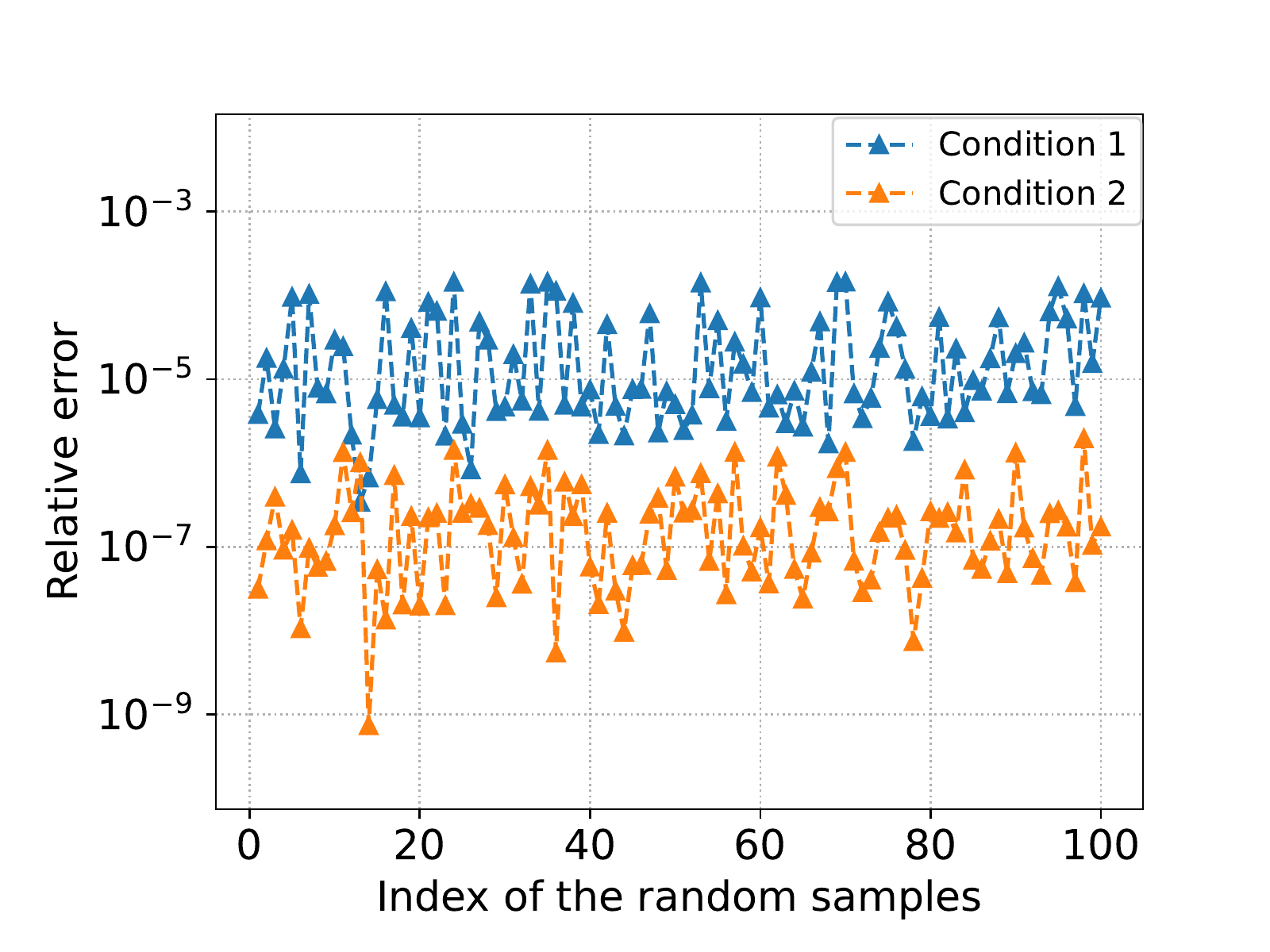}
	\caption{The relative error for 100 random samples by the SDD-VS method in two conditions: (1) $N_{S_1}=4, N_{S_2}=1, N_{F_1}=4, N_{F_2}=1$; (2) $N_{S_1}=6, N_{S_2}=1, N_{F_1}=4, N_{F_2}=1$.}
	\label{fig1}
\end{figure}

Based on the reduced model (\ref{eq-interf}), the interface solution by the SDD-VS method can be obtained by equation (\ref{eq-intso1}) directly. We implement it in two different conditions: (1) $N_{S_1}=4, N_{S_2}=1, N_{F_1}=4, N_{F_2}=1$; (2) $N_{S_1}=6, N_{S_2}=1, N_{F_1}=4, N_{F_2}=1$, and plot the relative errors in Figure \ref{fig1} to visualize the individual relative error of the first 100 samples from $N=10^4$ random samples. From the figure, we see that the method we proposed gives good approximation, and the relative error for each random sample becomes smaller when the number of the separated terms increases.

\begin{figure}[htbp]
	\centering
	\subfigure[$D_1$ ]{
		\includegraphics[width=3in, height=2.1in]{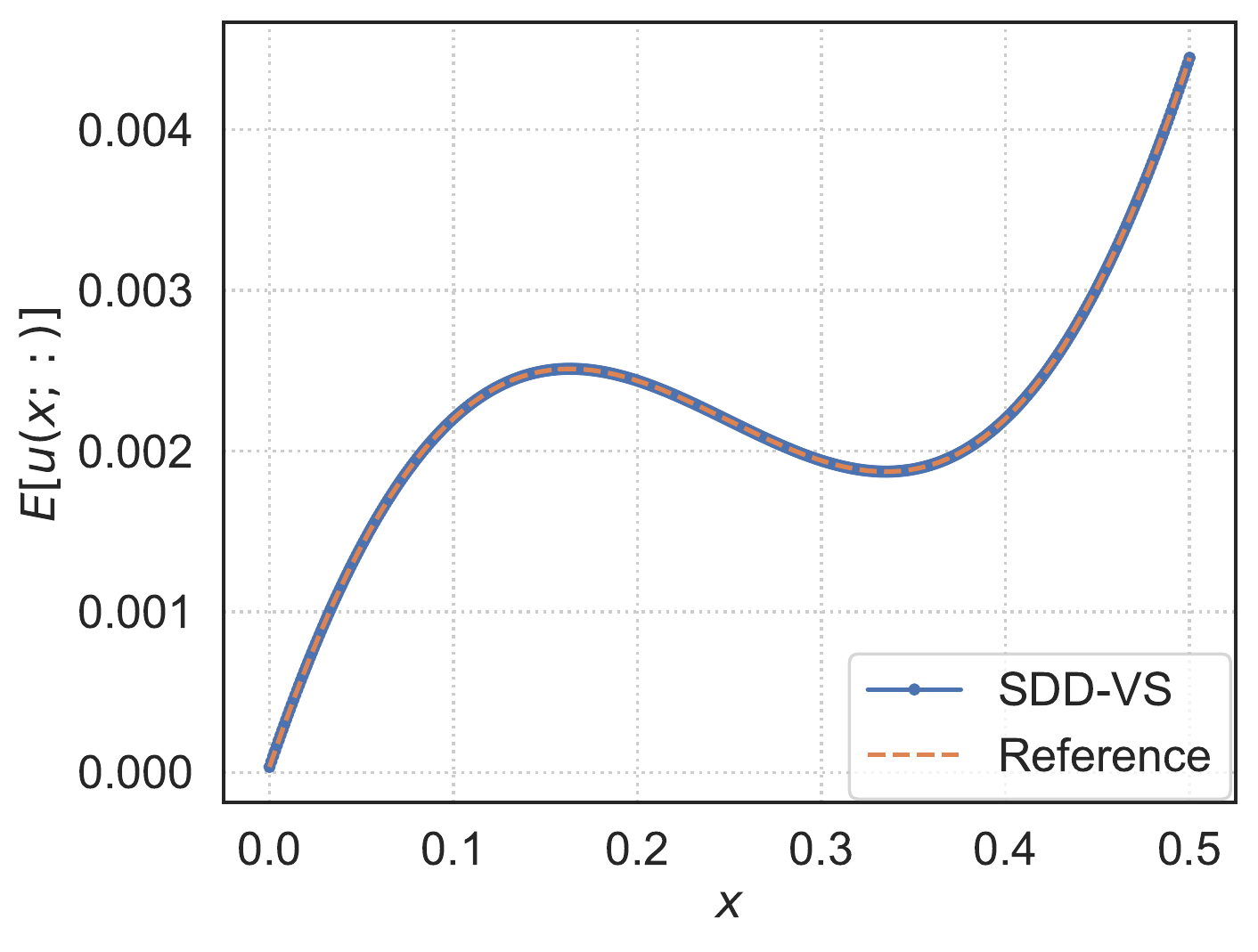}}
	\subfigure[$D_2$]{
		\includegraphics[width=3in, height=2.1in]{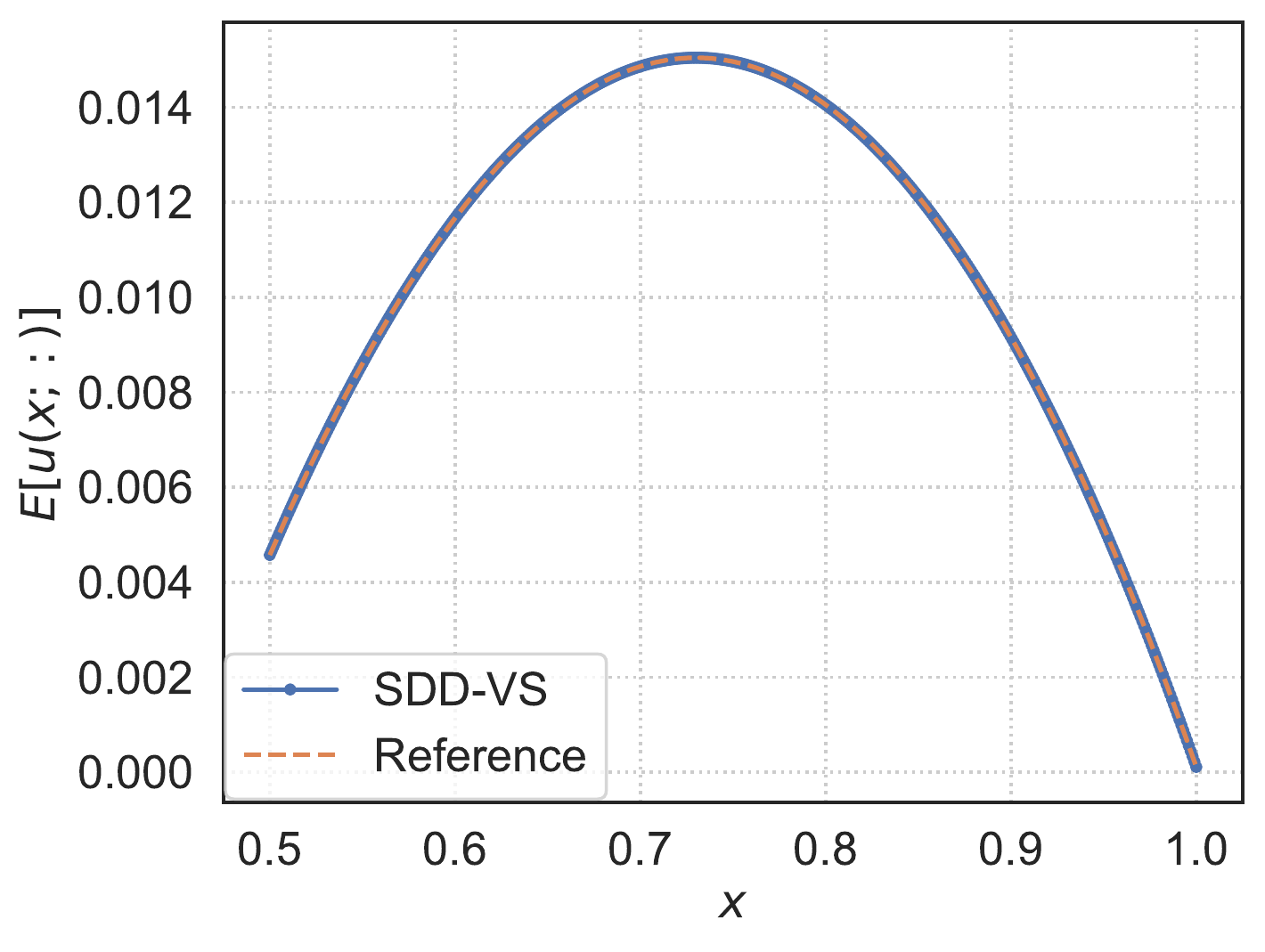}}
	\caption{Comparison of the mean solution for reference and the SDD-VS method with $N_{S_1}=4, N_{S_2}=1, N_{F_1}=4, N_{F_2}=1$ in subdomains $D_1$ and $D_2$.}
	\label{fig3}
\end{figure}

Figure \ref{fig3} shows the mean of the solution generated by our proposed SDD-VS method with $N_{S_1}=4, N_{S_2}=1, N_{F_1}=4, N_{F_2}=1$ and the reference solution defined on each subdomain, we find that our method provides good approximations.

\subsection{2D stochastic diffusion equation with three subdomains}
\label{num2}
In this example, the stochastic diffusion equation is posed on domain $D=[0,100]\times[0,100]$ consists of three subdomains, which is defined by
\begin{eqnarray*}
	\left\{
	\begin{aligned}
		-\text{div} \big(c(x;\bm{\xi}\big) \nabla u(x;\bm{\xi}))&=f(x;\bm{\xi}),  \ \forall ~x\in D, \\	
		u(x;\bm{\xi})&=g(x;\bm{\xi}),	 \ \forall ~x\in \partial D_1,\\
		\frac{\partial u(x;\bm{\xi})}{\partial n}&=h(x;\bm{\xi}), \ \forall ~x\in \partial D_2,
	\end{aligned}
	\right.
\end{eqnarray*}
where $\partial D_1=\{(0,0),(100,100)\}$, $\partial D_2=\partial D \backslash \partial D_1$, the source function $f(x;\bm{\xi})=0$, the Dirichlet boundary function $g(x;\bm{\xi})$ satisfied:  $g\big((0,0);\bm{\xi}\big)=20$,  $g\big((100,100);\bm{\xi}\big)=15$, the Neumann boundary function $h(x;\bm{\xi})=0$, and the random coefficient $c(x;\bm{\xi})$ is defined as
\begin{eqnarray*}
	c(x;\bm{\xi})=
	\left\{
	\begin{aligned}
		&\frac{80}{\mu}, &~\text{when} ~x\in D_1,\\
		&\frac{\bm{\xi}}{\mu}, &~\text{when} ~x\in D_2,\\
		&\frac{20}{\mu}, &~\text{when} ~x\in D_3,\\
	\end{aligned}
	\right.
\end{eqnarray*}
where $\mu=0.02$, $D_1=[0,100]\times\left[0,30\right]$, $D_2=[0,100]\times\left(30,70\right]$, $D_3$ is the remainder part of the domain $D$. We set the random variable $\bm{\xi}$ to be uniformly distributed in the interval $[1,4]$. Here, the finite element method calculates the reference solution with mesh size $h_x=h_y=1/100$.

For this numerical example, the affine expansion approximation of the stochastic stiffness matrix and the load vector defined on each subdomain has only one term, i.e. $m_{a_i}=m_{b_i}=1$, $i=1,2,3$, and then $N_{S_1}=N_{S_2}=N_{S_3}= N_{F_1}=N_{F_2}=N_{F_3}=1$. The affine expansion of $S(\bm{\xi})$ and $F(\bm{\xi})$ have a directly analytic expression with $m_S=2$, $m_F=1$ such as
\begin{eqnarray}
	\begin{aligned}
		S(\bm{\xi})&=\sum_{i=1}^{N_s}R_i^{T}S_i(\bm{\xi})R_i=\mathcal{\hat{X}}^1+\mathcal{\hat{X}}^2 \bm{\xi},\\
		F(\bm{\xi})&=\sum_{i=1}^{N_s}R_i^{T}F_i(\bm{\xi})={\hat{F}},
	\end{aligned}
\end{eqnarray}
where
$$\mathcal{\hat{X}}^1=\begin{bmatrix}
	A_{\Gamma_{12}\Gamma_{12}}^1-A_{\Gamma_{12} I}^1(A_{II}^1)^{-1}A_{I\Gamma_{12}}^1 &0\\ \quad
	0&A_{\Gamma_{23}\Gamma_{23}}^3-A_{\Gamma_{23} I}^3(A_{II}^3)^{-1}A_{I\Gamma_{23}}^3
\end{bmatrix},$$
$$\mathcal{\hat{X}}^2=\begin{bmatrix}
	A_{\Gamma_{12}\Gamma_{12}}^2-A_{\Gamma_{12} I}^2(A_{II}^2)^{-1}A_{I\Gamma_{12}}^2 &A_{\Gamma_{12}\Gamma_{23}}^2-A_{\Gamma_{12} I}^2(A_{II}^2)^{-1}A_{I\Gamma_{23}}^2\\ \\
	
	A_{\Gamma_{23}\Gamma_{12}}^2-A_{\Gamma_{23} I}^2(A_{II}^2)^{-1}A_{I\Gamma_{12}}^2&A_{\Gamma_{23}\Gamma_{23}}^2-A_{\Gamma_{23} I}^2(A_{II}^2)^{-1}A_{I\Gamma_{23}}^2
\end{bmatrix},$$
and
$${\hat{F}}=
\begin{Bmatrix}
	\bm{f}_{\Gamma_{12}}^1-A_{\Gamma_{12} I}^1(A_{II}^1)^{-1}\bm{f}_{I}^1+\bm{f}_{\Gamma_{12}}^2-A_{\Gamma_{12} I}^2(A_{II}^2)^{-1}\bm{f}_{I}^2\\ \\
	\bm{f}_{\Gamma_{23}}^2-A_{\Gamma_{23} I}^2(A_{II}^2)^{-1}\bm{f}_{I}^2+\bm{f}_{\Gamma_{23}}^3-A_{\Gamma_{23} I}^3(A_{II}^3)^{-1}\bm{f}_{I}^3
\end{Bmatrix},$$
are independent
of random variable $\bm{\xi}$, and their computation is once.
We choose $|\Xi|=20$ samples to get the efficient surrogate model of $\bm{u}_{\Gamma}(\bm{\xi})$ by Algorithm \ref{algorithm-vs-sSAS} based on the following stochastic algebraic system $$(\mathcal{\hat{X}}^1+\mathcal{\hat{X}}^2 \bm{\xi}) \bm{u}_{\Gamma}(\bm{\xi})={\hat{F}}.$$

\begin{figure}[htbp]
	\centering
	\includegraphics[width=3.8in, height=2.7in]{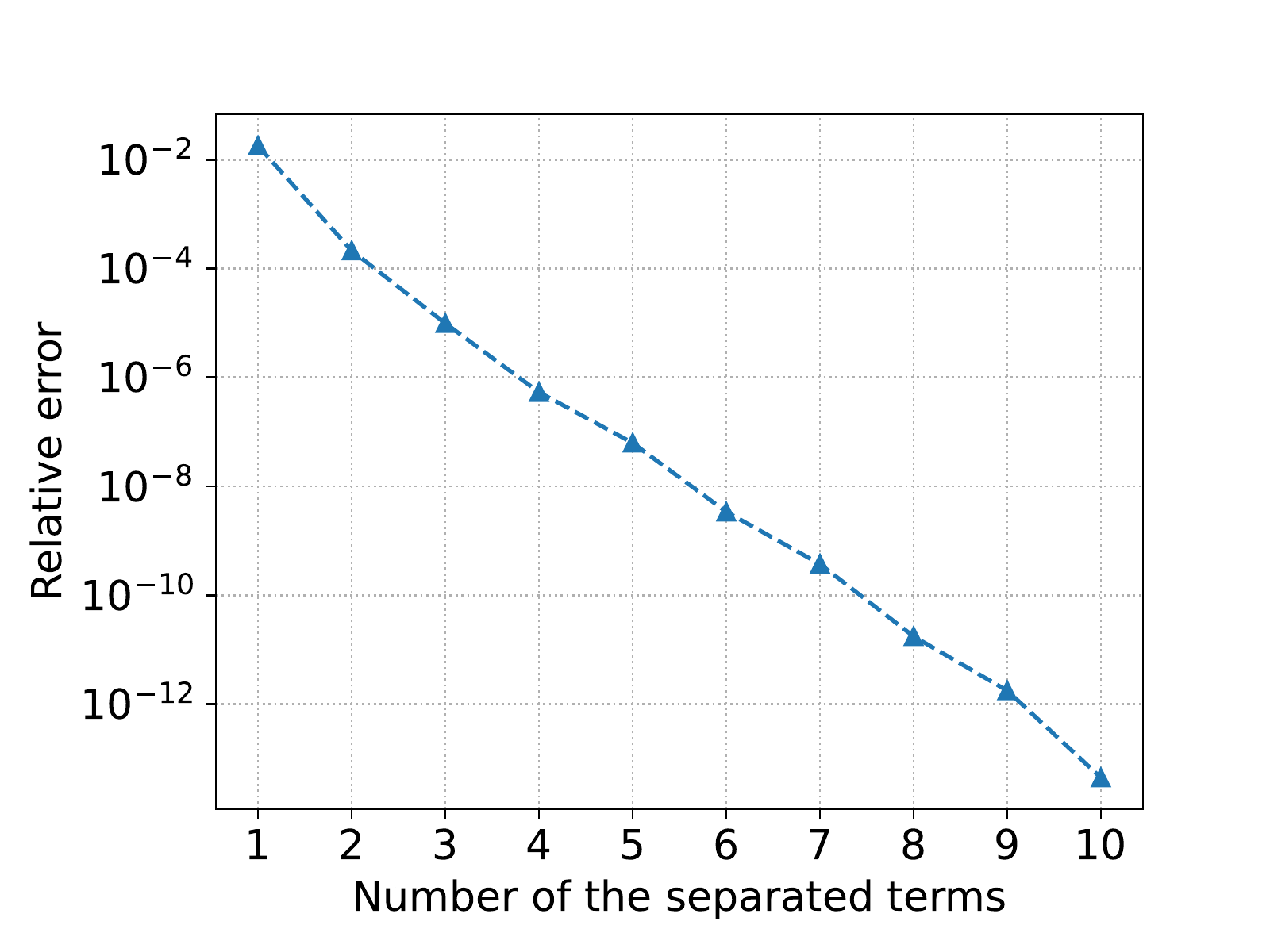}
	\caption{Comparison of the relative mean error corresponding to the different number of the separated terms $M$.}
	\label{fig2.1}
\end{figure}

Firstly, we focus on the results of the stochastic interface problem solved by the SDD-VS method.
In Figure \ref{fig2.1}, we depict the relative mean error of the stochastic interface problem by the SDD-VS method versus a different number of the separated terms $M$, where the relative mean error $\epsilon$ is calculated with $N=10^4$ test samples. According to the figure, we can see that as the number of the separated terms $M$ increases, the approximation becomes more accurate.

\begin{figure}[htbp]
	\centering
	\includegraphics[width=3.8in, height=2.7in]{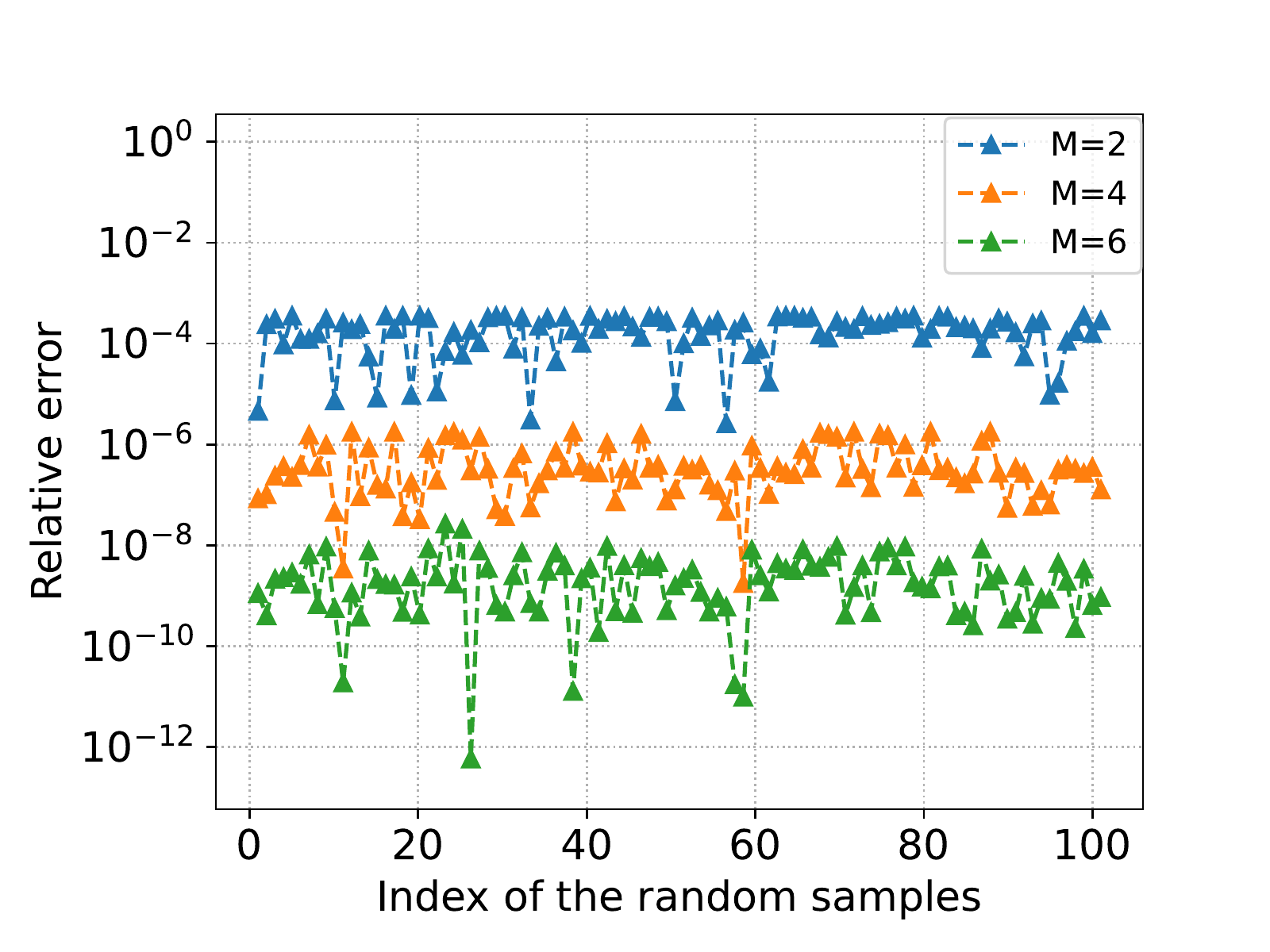}
	\caption{The relative error for 100 samples versus the number of the separated terms being $M=2$, $M=4$, $M=6$.}
	\label{fig2.2}
\end{figure}

To visualize the individual relative error, we choose the first $100$ samples from $N=10^4$ random samples and plot the relative error with the number of the separated terms being $M=2, M=4, M=6$ in Figure \ref{fig2.2}. This shows that the relative error for each individual sample becomes smaller when the number of the separated terms $M$ increases.

\begin{figure}[htbp]
	\centering
	\subfigure[$\Gamma_1$ ]{
		\includegraphics[width=3in, height=2.1in]{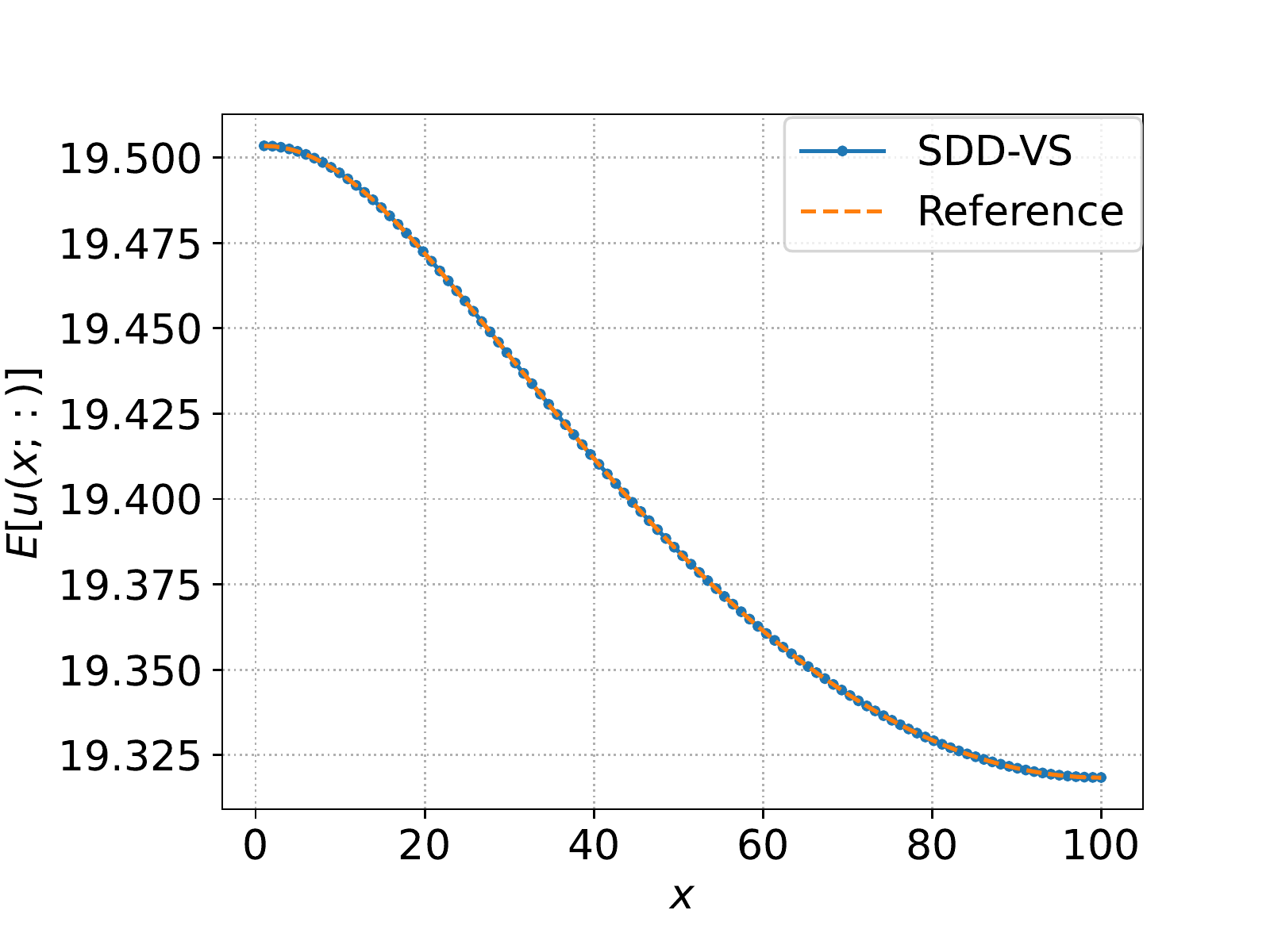}}
	\subfigure[$\Gamma_2$]{
		\includegraphics[width=2.8in, height=1.9in]{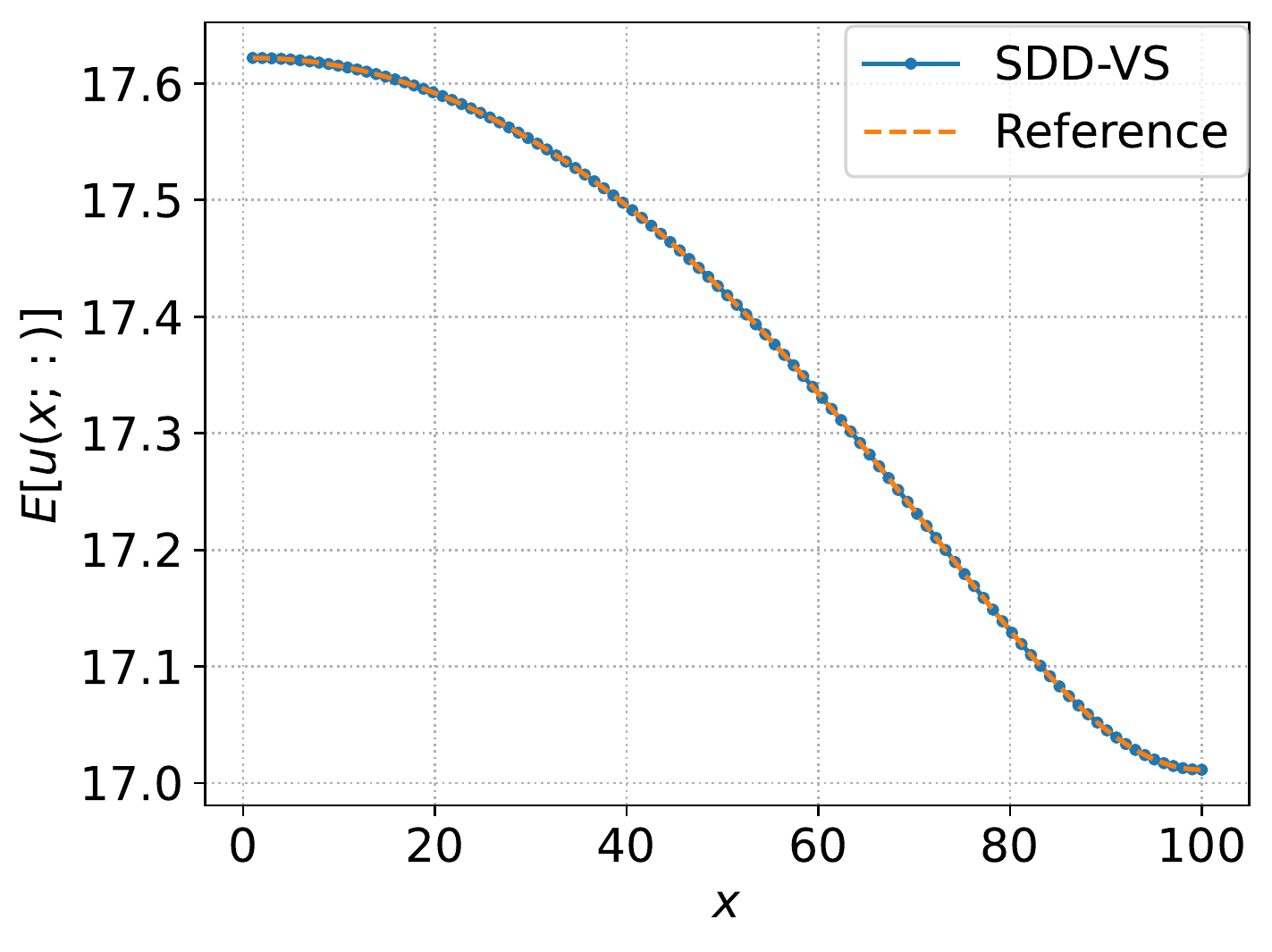}}
	\caption{Comparison of the mean solution for the reference and the SDD-VS method in two interfaces with the number of the separated terms $M=3$.}
	\label{fig2.3}
\end{figure}

Figure \ref{fig2.3} shows the mean of the solution generated by the SDD-VS method and the reference solution in two interfaces $[0,100]\times30$ and $[0,100]\times70$ with the number of the separated terms $M=3$ based on $10^4$ random samples, as we can see, both interfaces yield an accurate mean solution.

\begin{figure}[htbp]
	\centering
	\subfigure[]{
		\includegraphics[width=3in, height=2.1in]{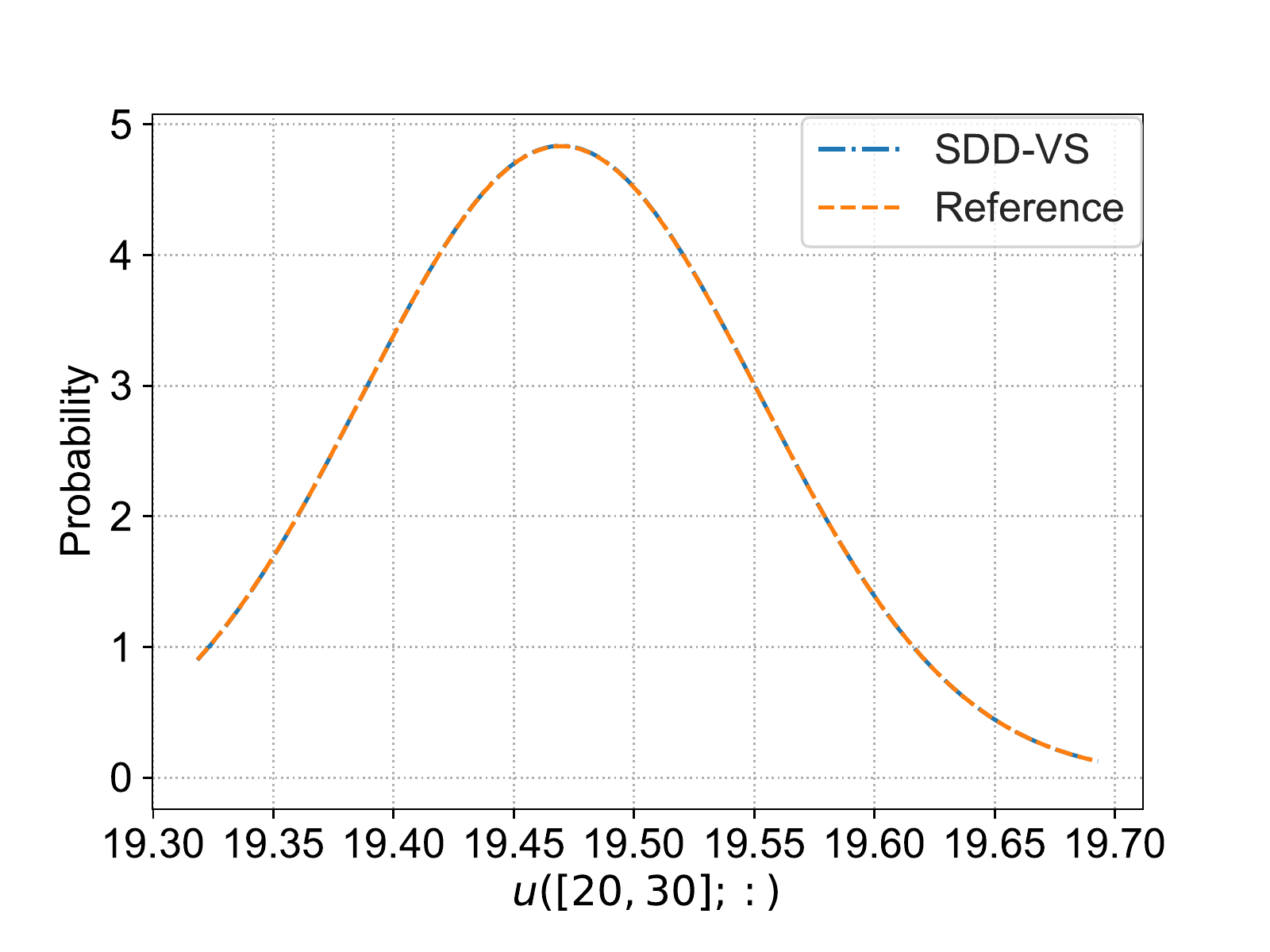}}
	\subfigure[]{
		\includegraphics[width=3in, height=2.1in]{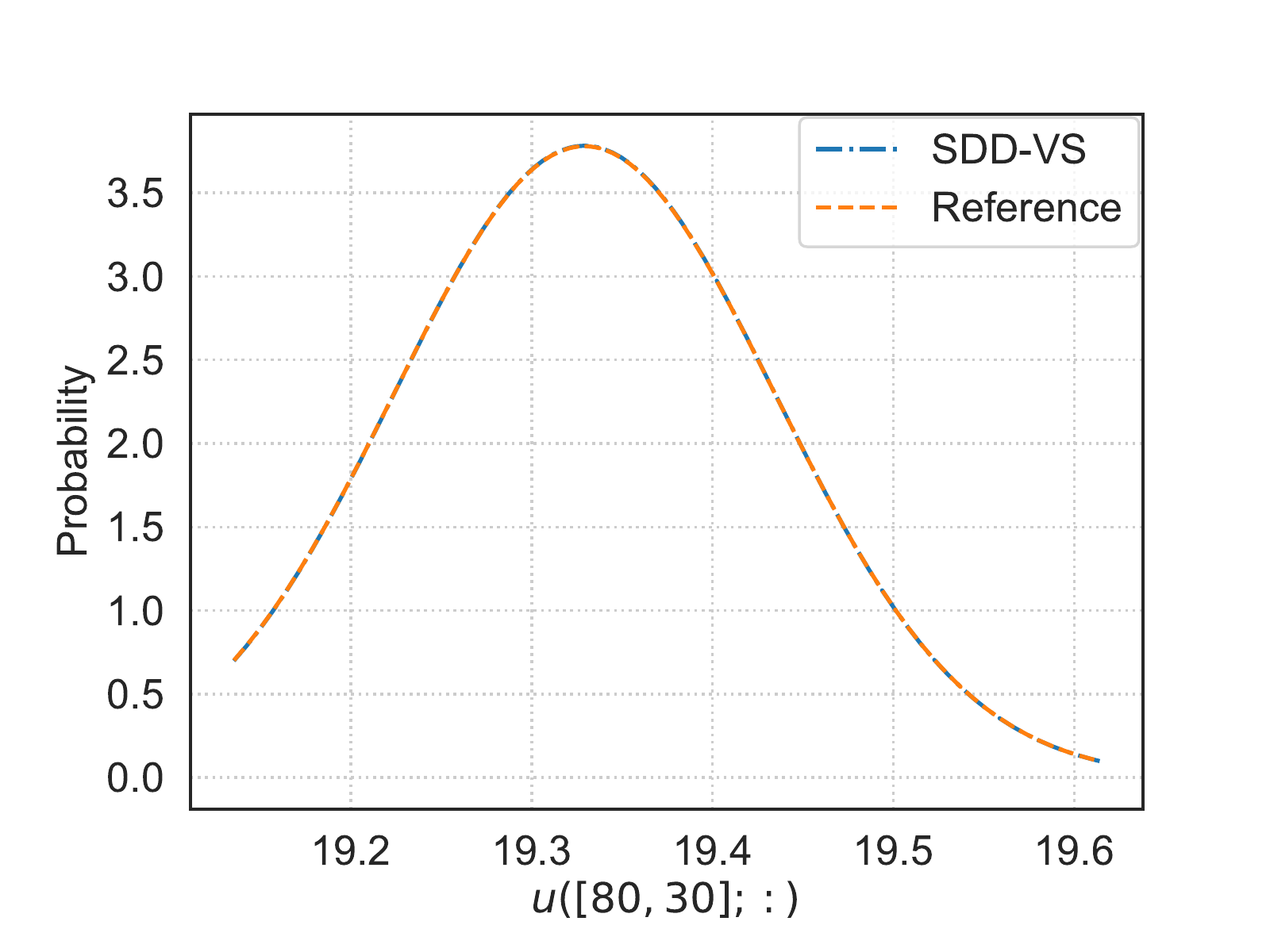}}
	\\
	\subfigure[]{
		\includegraphics[width=3in, height=2.1in]{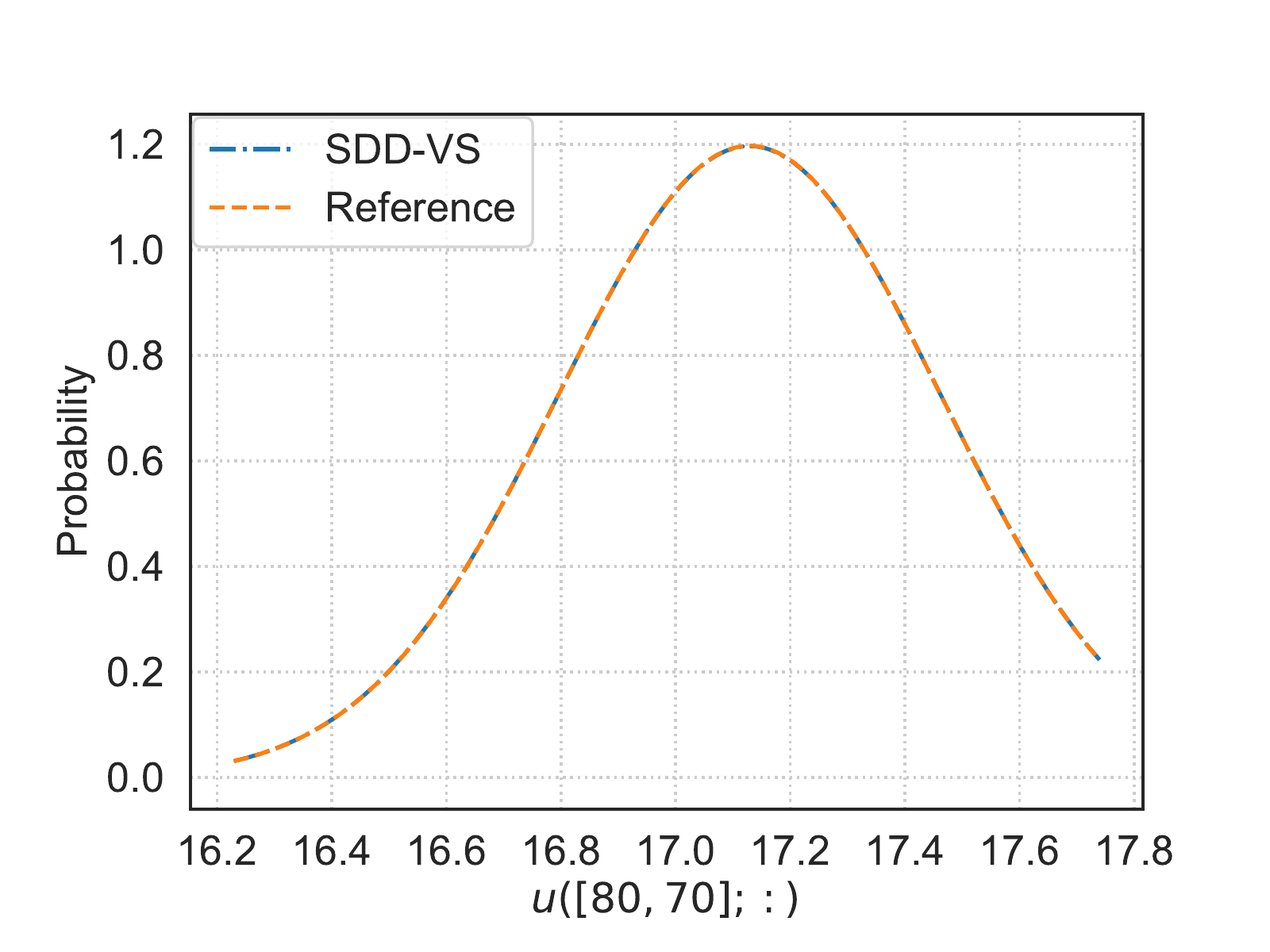}}
	\subfigure[]{
		\includegraphics[width=3in, height=2.1in]{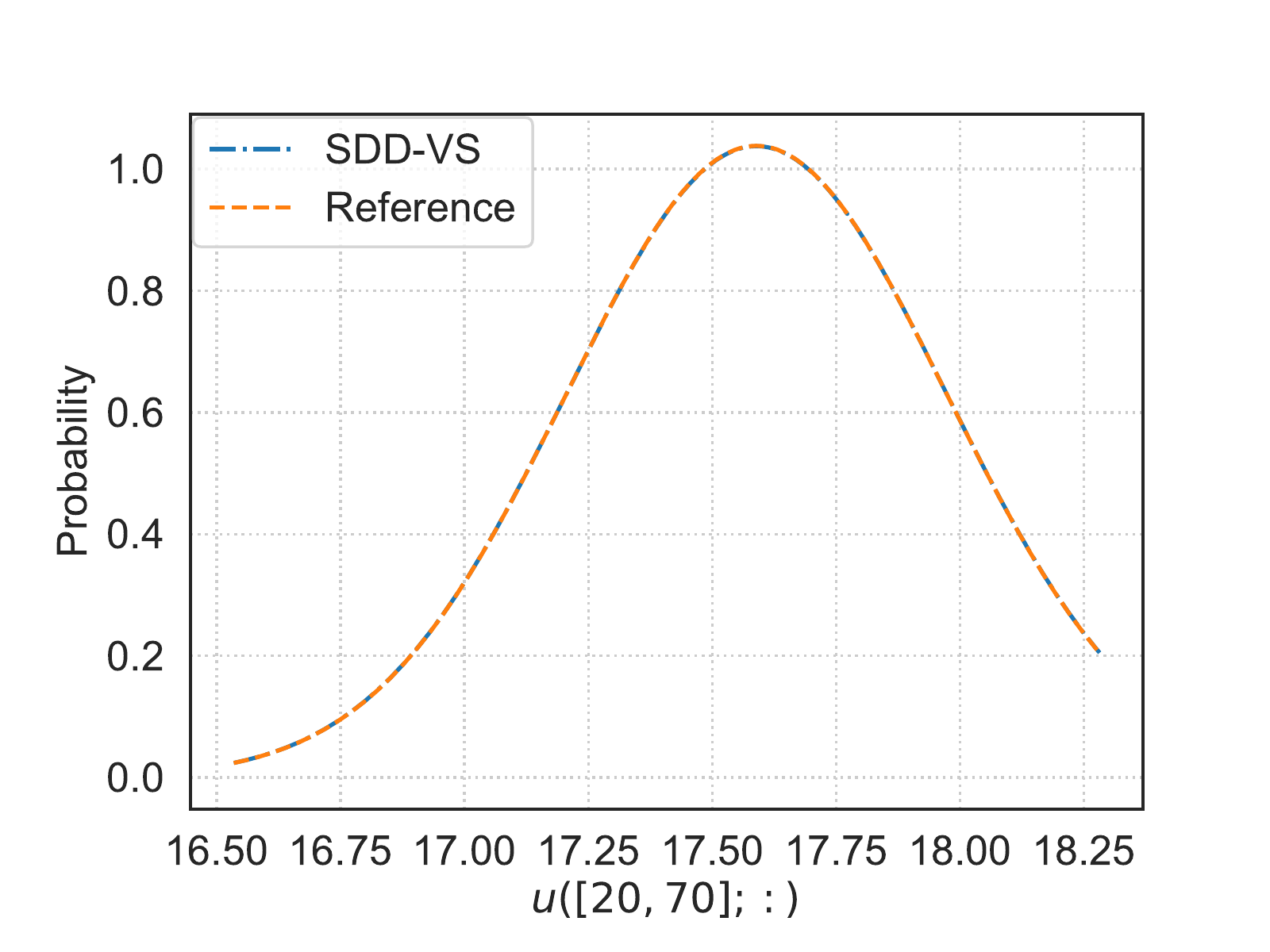}}
	\caption{Comparison of the probability density of $u(x_0;\bm{\xi})$ for reference and the SDD-VS method with the number of the separated terms $M=3$.}
	\label{fig2.4}
\end{figure}

The probability density estimates of the reference and the SDD-VS method at a single measurement location in interfaces are shown in Figure \ref{fig2.4}. From the figure, we find that the SDD-VS method can give a good approximation for the reference probability density.

Following (\ref{eq-interior}), the stochastic interior unknowns $\bm{u}_{I}^i(\bm{\xi})$ belong to each subdomain $D_i$ can be obtained by solving the following stochastic linear equation
\begin{eqnarray*}
	A_{II}^{i}(\bm{\xi})\bm{u}_I^i(\bm{\xi})=\bm{f}_I^i(\bm{\xi})-A_{I\Gamma}^{i}(\bm{\xi})\bm{u}_{\Gamma}^i(\bm{\xi}).
\end{eqnarray*}

\begin{figure}[htbp]
	\centering
	\subfigure[$D_1$ ]{
		\includegraphics[width=1.6in, height=1.5in]{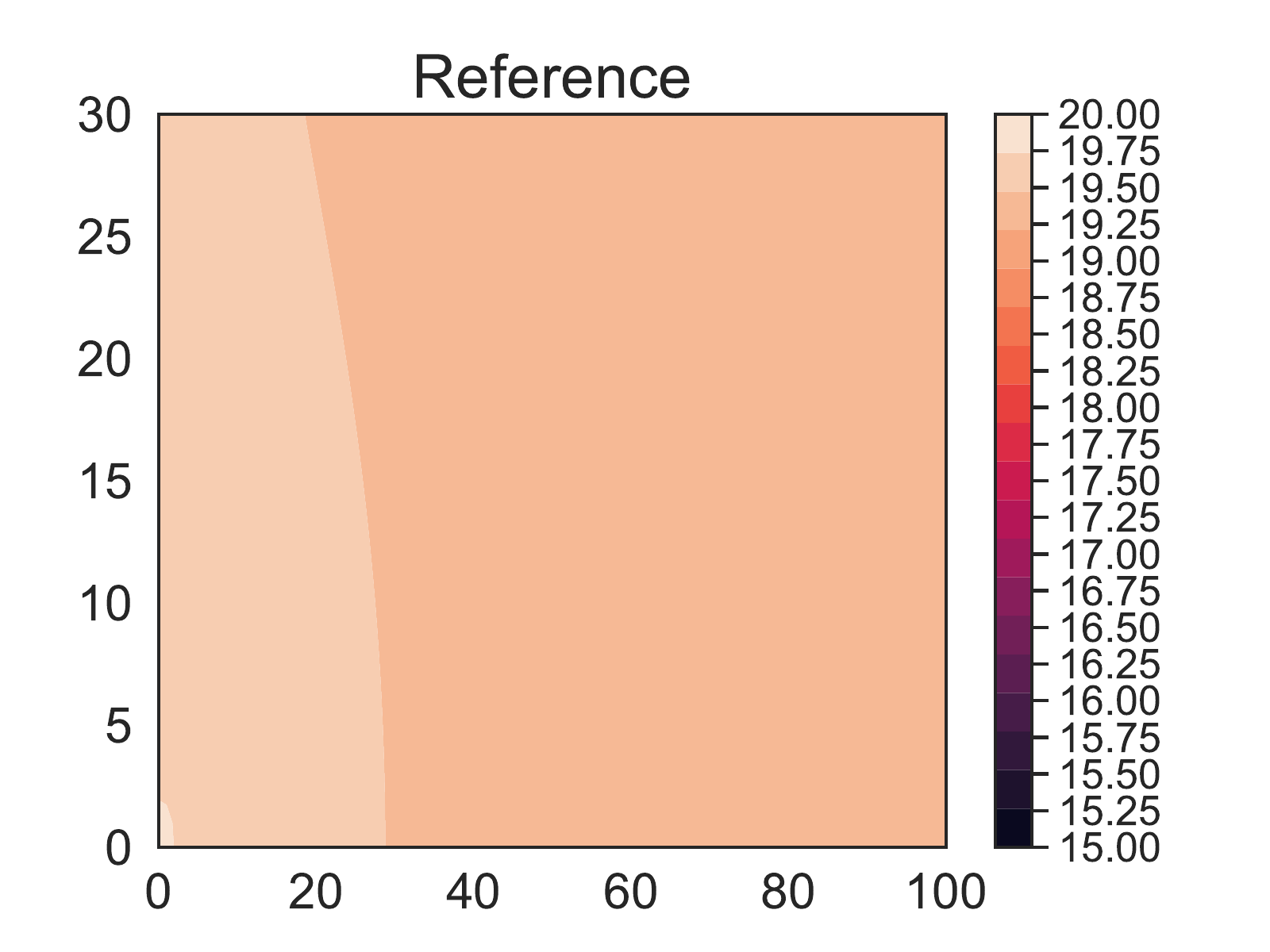}}
	\subfigure[$D_2$]{
		\includegraphics[width=1.6in, height=1.5in]{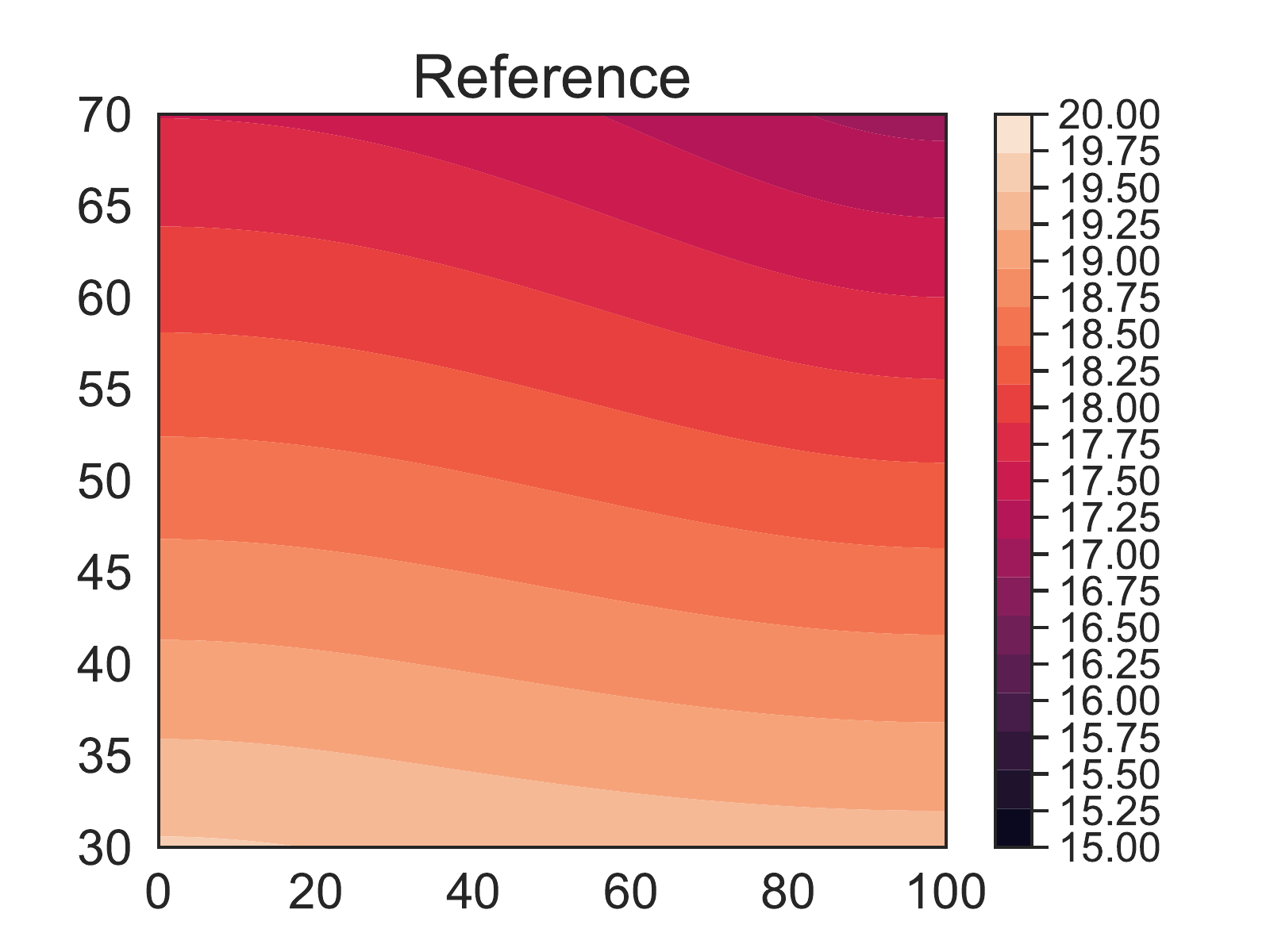}}
	\subfigure[$D_3$]{
		\includegraphics[width=1.6in, height=1.5in]{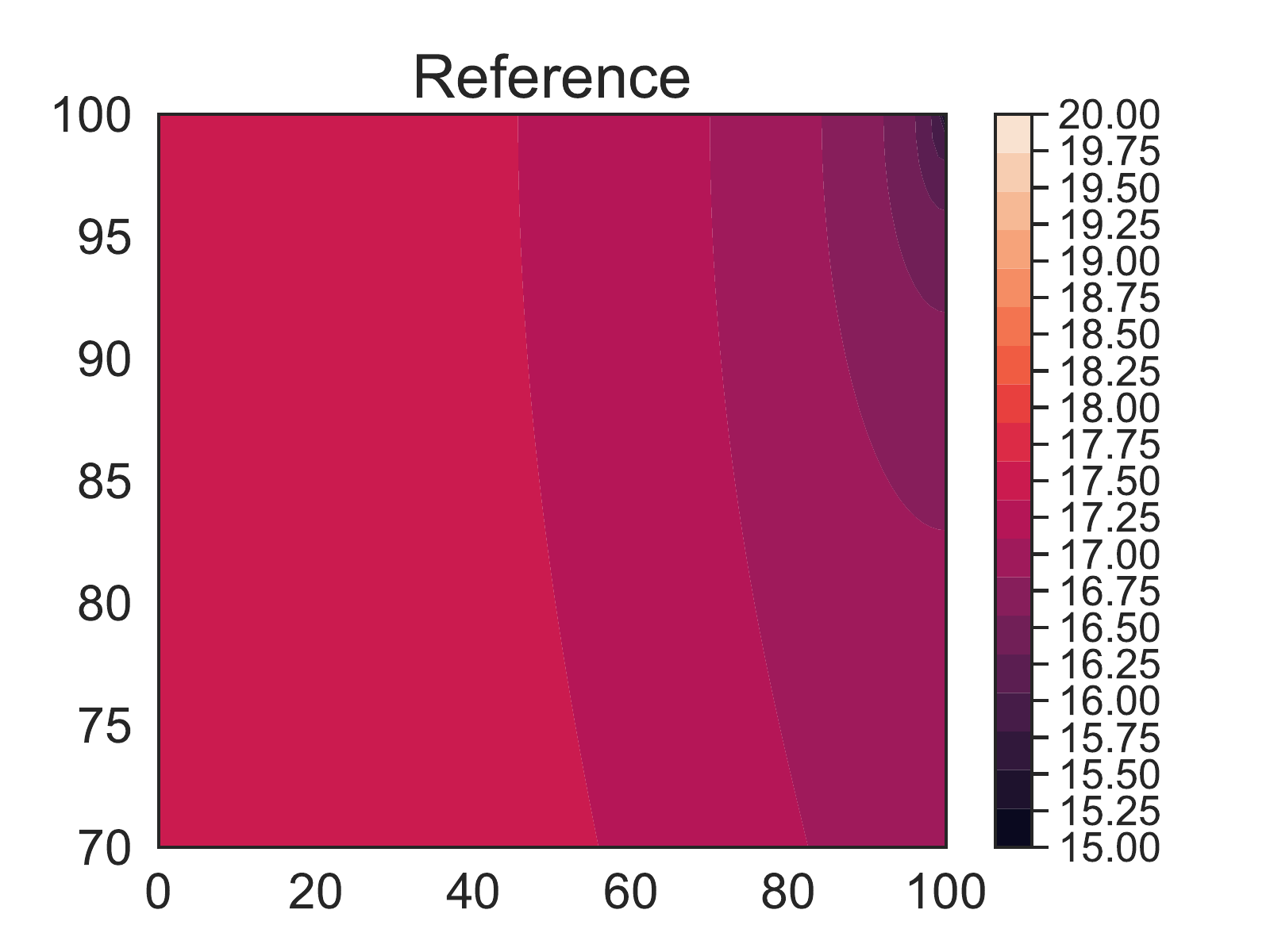}}\\
	\subfigure[$D_1$]{
		\includegraphics[width=1.6in, height=1.5in]{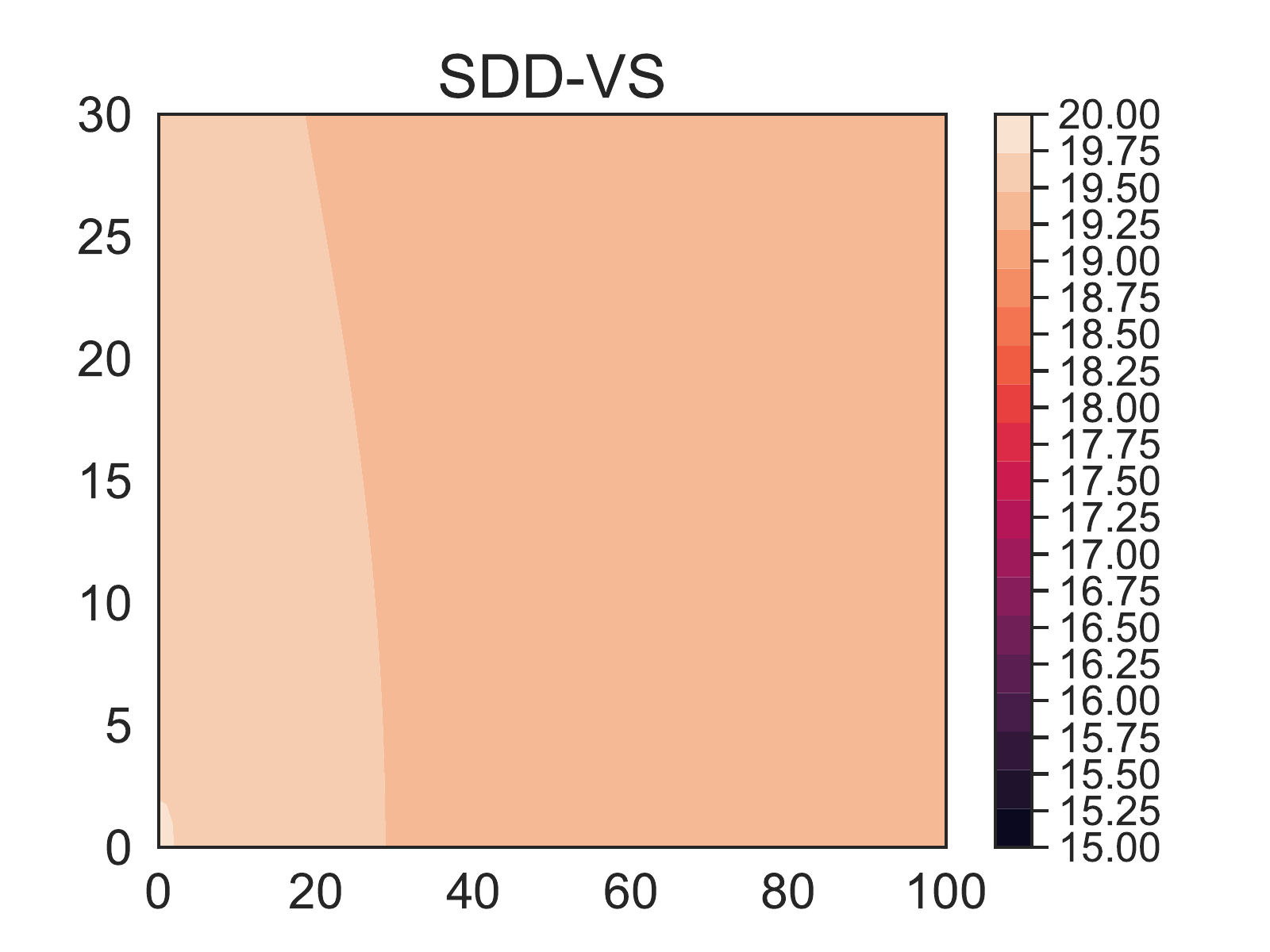}}
	\subfigure[$D_2$]{
		\includegraphics[width=1.6in, height=1.5in]{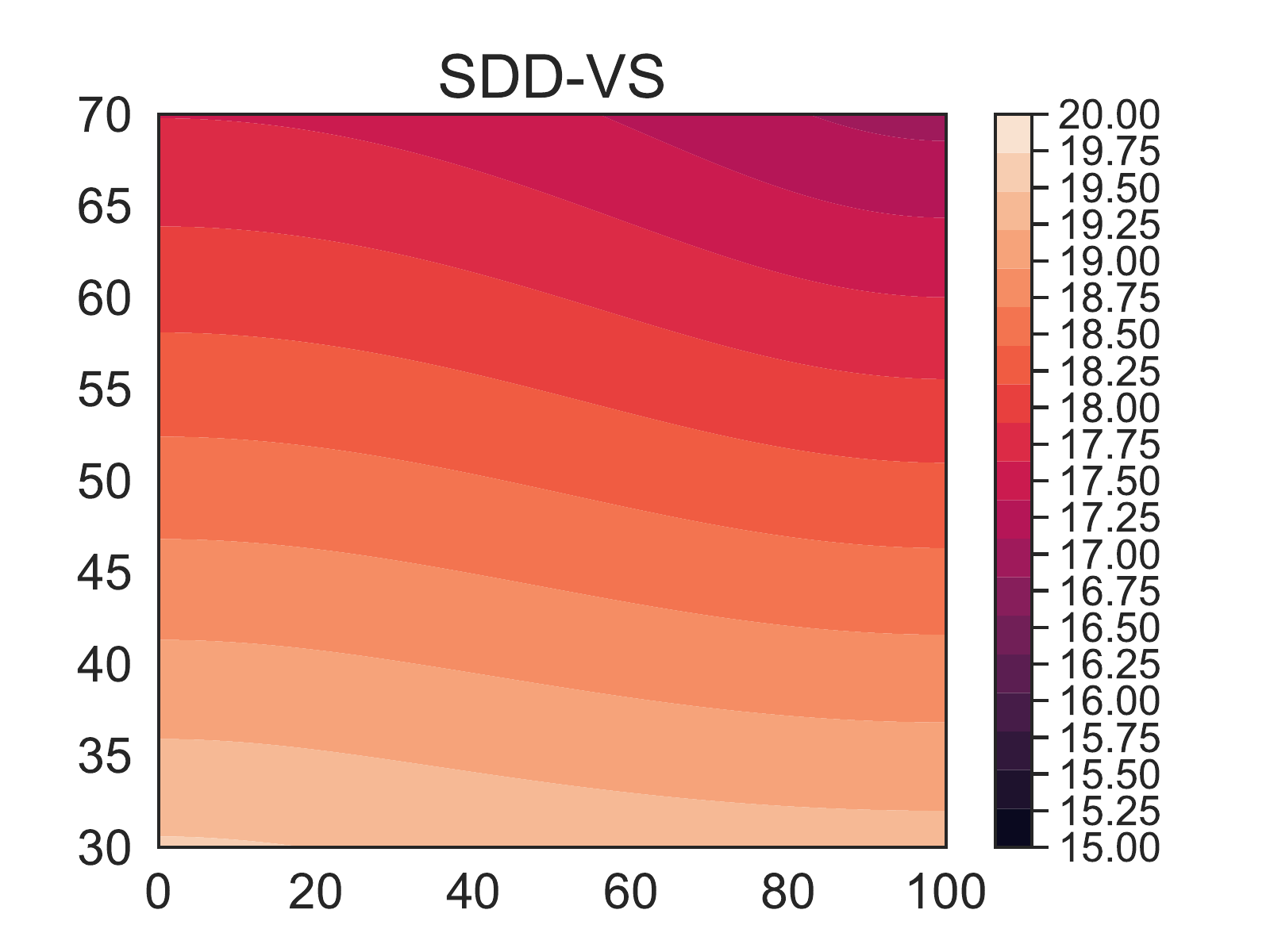}}
	\subfigure[$D_3$]{
		\includegraphics[width=1.4in, height=1.5in]{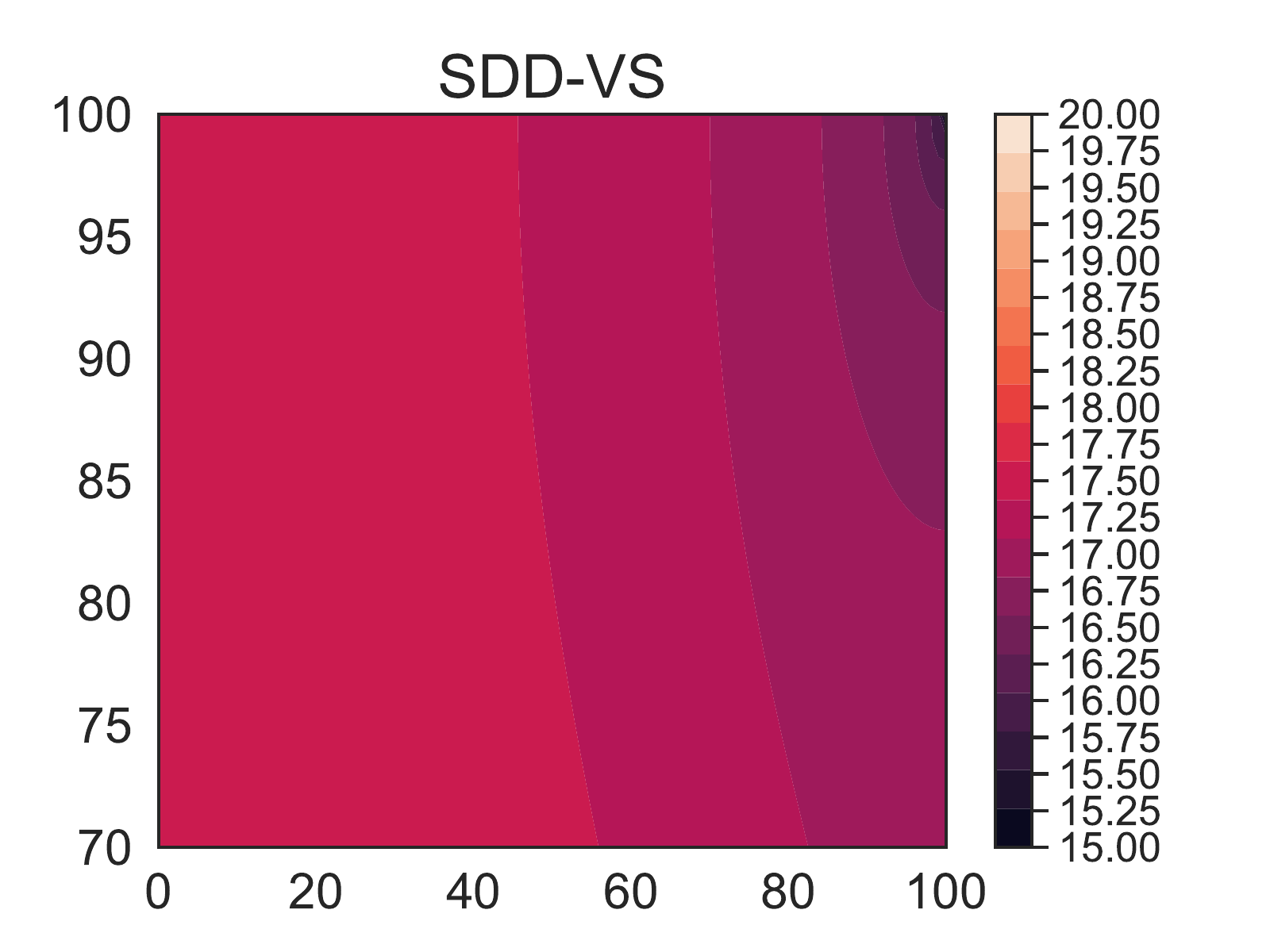}}
	\caption{Comparison of the mean solution on three subdomains for reference and the SDD-VS method with the number of the separated terms $M=3$.}
	\label{fig2.5}
\end{figure}

Indeed, for this type of SPDE, the subproblems can be solved more efficiently. In subdomain $D_i$, once we get the reduced model representation of the stochastic interface unknowns $\bm{u}_{\Gamma}^{i}(\bm{\xi})$, the reduced model representation of the interior unknowns $\bm{u}_{I}^i(\bm{\xi})$ can be obtained directly, which can be written as
\begin{eqnarray}
	\label{eq-sinter1}
	\begin{aligned}
		\bm{u}_I^i(\bm{\xi})&=(A_{II}^{i})^{-1}(\bm{f}_I^i(\bm{\xi})-A_{I\Gamma}^{i}\bm{u}_{\Gamma}^i(\bm{\xi}))\\
		&\approx(A_{II}^{i})^{-1}\sum_{j=1}^{m_{b_i}}q^{ij}(\bm{\xi})\bm{f}_I^{ij}-(A_{II}^{i})^{-1}A_{I\Gamma}^{i}\sum_{j=1}^{M}\zeta_j(\bm{\xi})\bm{c}_{j}^{i}\\
		&=\sum_{j=1}^{m_{b_i}}q^{ij}(\bm{\xi})\hat{\bm{f}}_I^{ij}-\sum_{j=1}^{M}\zeta_j(\bm{\xi})\hat{\bm{c}}_{j}^{i}, ~i=1,3,\\
	\end{aligned}
\end{eqnarray}
\begin{eqnarray}
	\label{eq-sinter2}
	\begin{aligned}
		\bm{u}_I^i(\bm{\xi})&=(A_{II}^{i})^{-1}\Big(\frac{\bm{f}_I^i(\bm{\xi})}{\bm{\xi}}-A_{I\Gamma}^{i}\bm{u}_{\Gamma}^i(\bm{\xi})\Big)\\
		&\approx(A_{II}^{i})^{-1}\sum_{j=1}^{m_{b_i}}\frac{q^{ij}(\bm{\xi})}{\bm{\xi}}\bm{f}_I^{ij}-(A_{II}^{i})^{-1}A_{I\Gamma}^{i}\sum_{j=1}^{M}\zeta_j(\bm{\xi})\bm{c}_{j}^{i}\\
		&=\sum_{j=1}^{m_{b_i}}\frac{q^{ij}(\bm{\xi})}{\bm{\xi}}\hat{\bm{f}}_I^{ij}-\sum_{j=1}^{M}\zeta_j(\bm{\xi})\hat{\bm{c}}_{j}^{i}, ~i=2,
	\end{aligned}
\end{eqnarray}
where the second equation following from $A_{II}^{i}(\bm{\xi})=A_{II}^{i}\bm{\xi}$, $A_{I\Gamma}^{i}(\bm{\xi})=A_{I\Gamma}^{i}\bm{\xi}$, when $i=2$, for each $i$, $\bm{c}_{j}^{i}$ is the corresponding vector of $\bm{c}_{j}$ restricted on the interface of $D_i$, $\hat{\bm{f}}_I^{ij}=(A_{II}^{i})^{-1}\bm{f}_I^{ij}$, $\hat{\bm{c}}_{j}^{i}=(A_{II}^{i})^{-1}A_{I\Gamma}^{i}\bm{c}_{j}^{i}$. In Figure (\ref{fig2.5}), we plot the mean of the reference solution and the solution generated by (\ref{eq-sinter1}-\ref{eq-sinter2}) in three subdomains with the number of the separated terms $M=3$ for the stochastic interface unknowns, the first row is the reference solution, and the second row is the solution generated by our proposed SDD-VS method. We find that our method works well for this problem.

Finally, since the VS method proposed in \cite{Li2017a} is applicable to solve SPDEs directly, we compare the accuracy and the computational cost of the SDD-VS method against the VS method. To this end, $|\Xi|=20$ samples are selected to construct the surrogate model for the VS method.
In Table \ref{table3_1}, we list the relative mean error, CPU time (offline CPU time $\mathcal{T}_{\text{off}}$, online CPU time $\mathcal{T}_{\text{on}}$, total CPU time $\mathcal{T}_{\text{tot}}$ and average online CPU time $T_{\text{on}}$) for SDD-VS method with the number of the separated terms $M=10$ for the interface unknowns, VS method with the number of the separated terms $M_{VS}=3,6,12$ and the reference method.
Note that the online CPU time for the SDD-VS method can become smaller if we use the parallel strategy. From the table, we conclude that
(1) as the number of separated terms $M_{VS}$ increases, the CPU times
needed for the VS method increase steadily;
(2) the magnitude of average online CPU time by the SDD-VS method and VS method are much smaller than that of the reference method;
(3)
the SDD-VS method achieves much better approximation than the VS method with similar online computational cost, besides, the SDD-VS method uses much less computational cost than the VS method with similar approximation accuracy.
In summary, the proposed SDD-VS method for this problem renders a more robust and accurate approximation than the VS method applied here directly.
\begin{table}[hbtp]
	\centering
	\caption{Comparison of average relative errors and the CPU times for SDD-VS, VS and the reference based on $10^4$ parameter samples.}
	\vspace*{2pt}
	\begin{tabular}{c|c|c|c|c|c}
		\Xhline{1pt}
		\multirow{1}{*}{\parbox[t]{1.5cm}{\centering }} &\multicolumn{1}{c|}{} &
		\multicolumn{3}{c|}{}\\
		&
		\multicolumn{1}{c|}{\centering SDD-VS }
		&
		\multicolumn{3}{c|}{\centering VS}
		&
		\multicolumn{1}{c}{\centering Reference}\\
		\cline{3-5} \cline{3-5}
		& \hspace*{1cm}& \hspace*{1cm} & \hspace*{1cm} & \hspace*{1cm}& \hspace*{1cm}   \\[-9pt]
		& & $M_{VS}=3$ & $M_{VS}=6$ & $M_{VS}=12$ &  \\
		\hline
		$\varepsilon$  &  $2.40\times10^{-13}$&  $1.37\times10^{-5}$ & $8.14\times10^{-9}$& $1.37\times10^{-13}$   & $\setminus$ \\
		\hline
		$\mathcal{T}_{\text{off}}$&$10.39s$  &  $222.14 s$ & $223.23s$& $225.95s$  & $\setminus$ \\
		\hline
		$\mathcal{T}_{\text{on}}$&$3.98s$  &  $9.53s$ & $22.65s$& $59.67s$ & $\setminus$ \\
		\hline
		$\mathcal{T}_{\text{tot}}$&$14.37s$  &  $231.67 s$ & $245.88 s$ & $285.62s$  & $94507.32s$\\
		\hline
		$T_{\text{on}}$ &$3.98\times10^{-4}s$ &  $9.53\times10^{-4}s$ & $2.27\times10^{-3} s$& $5.97\times10^{-3}s$  & $9.45s$\\
		\Xhline{1pt}
	\end{tabular}
	\label{table3_1}
\end{table}

\subsection{2D stochastic convection-diffusion equation with high-dimensional random inputs}
\label{num3}
In the final example, we consider the stochastic convection-diffusion equation defined on domain $D=[0,1]\times[0,1]$ with high-dimensional random variables, which is defined by
\begin{eqnarray*}
	\left\{
	\begin{aligned}
		-\text{div} \big(c(x,y;\bm{\xi}\big) \nabla u(x,y;\bm{\xi}))+d(x,y;\bm{\xi}) \cdot\nabla u(x,y;\bm{\xi}) &=f(x.y;\bm{\xi}),  \ \forall ~(x,y)\in D, \\	
		u(x,y;\bm{\xi})&=g(x,y;\bm{\xi}),	 \ \forall ~(x,y)\in \partial D,
	\end{aligned}
	\right.
\end{eqnarray*}
where the boundary condition function $g(x,y;\bm{\xi})=0$, the random diffusivity $c(x,y;\bm{\xi})$ and velocity $d(x,y;\bm{\xi})$ are defined as
\begin{eqnarray*}
	c(x,y;\bm{\xi})=\bm{\xi}_{33}+y\bm{\xi}_{34}+3, ~~d(x,y;\bm{\xi})=1.
\end{eqnarray*}

For this problem, we divide $D$ into two subdomains $D_1=[0,0.5]\times[0,1]$ and $D_2=\left(0.5,1\right]\times [0,1]$.
{In each subdomain $D_i$, the source function $f(x,y;\bm{\xi})$ is taken as a random field, which is characterized by a two-point exponential covariance function cov$[f]$, i.e.,
	\begin{eqnarray}
		\label{eq-cov}
		\text{cov}[f](x_1,y_1;x_2,y_2)=\sigma^2\text{exp}\big(-\frac{|x_1-x_2|^2}{2l_x^2}-\frac{|y_1-y_2|^2}{2l_y^2}\big),
	\end{eqnarray}
	where $(x_i,y_i),i=1,2$ is the spatial coordinate, the variance $\sigma=0.1$, correlation length $l_y=0.5$ and we take different correlation lengths in the $x$ direction for different subdomains, specifically, $l_x=0.5$ in subdomain $D_1$ while $l_x=0.05$ in subdomain $D_2$.}
Then the random source term $f(x,y;\bm{\xi})$ is obtained by truncating a Karhunen-Lo$\grave{e}$ve expansion, i.e.,
\begin{eqnarray*}
	f(x,y;\bm{\xi}):=E[f]+\sum_{i=1}^{16}\sqrt{\gamma_i}b_i(x,y)\bm{\xi}_i,
\end{eqnarray*}
where $E[f]=1$,
the random vector $\bm{\xi}:=({\xi}_1,{\xi}_2,\cdots,{\xi}_{34})\in\mathbb{R}^{34} $, we assumed that ${\xi}_i,i=1,2,\cdots,34$ are i.i.d uniform random variables range $[-1,1]$. The reference solution is obtained by the finite element method with mesh size $h_x=h_y=1/60$.

We focus on the numerical result of the stochastic interface problem solved by the SDD-VS method with  $N_{S_1}=N_{S_2}=20$, $N_{F_1}=N_{F_2}=80$. And we choose $|\Xi|=120$ samples for training in Algorithm \ref{algorithm-vs-sSAS}.

\begin{figure}[htbp]
	\centering
	\includegraphics[width=3.8in, height=2.7in]{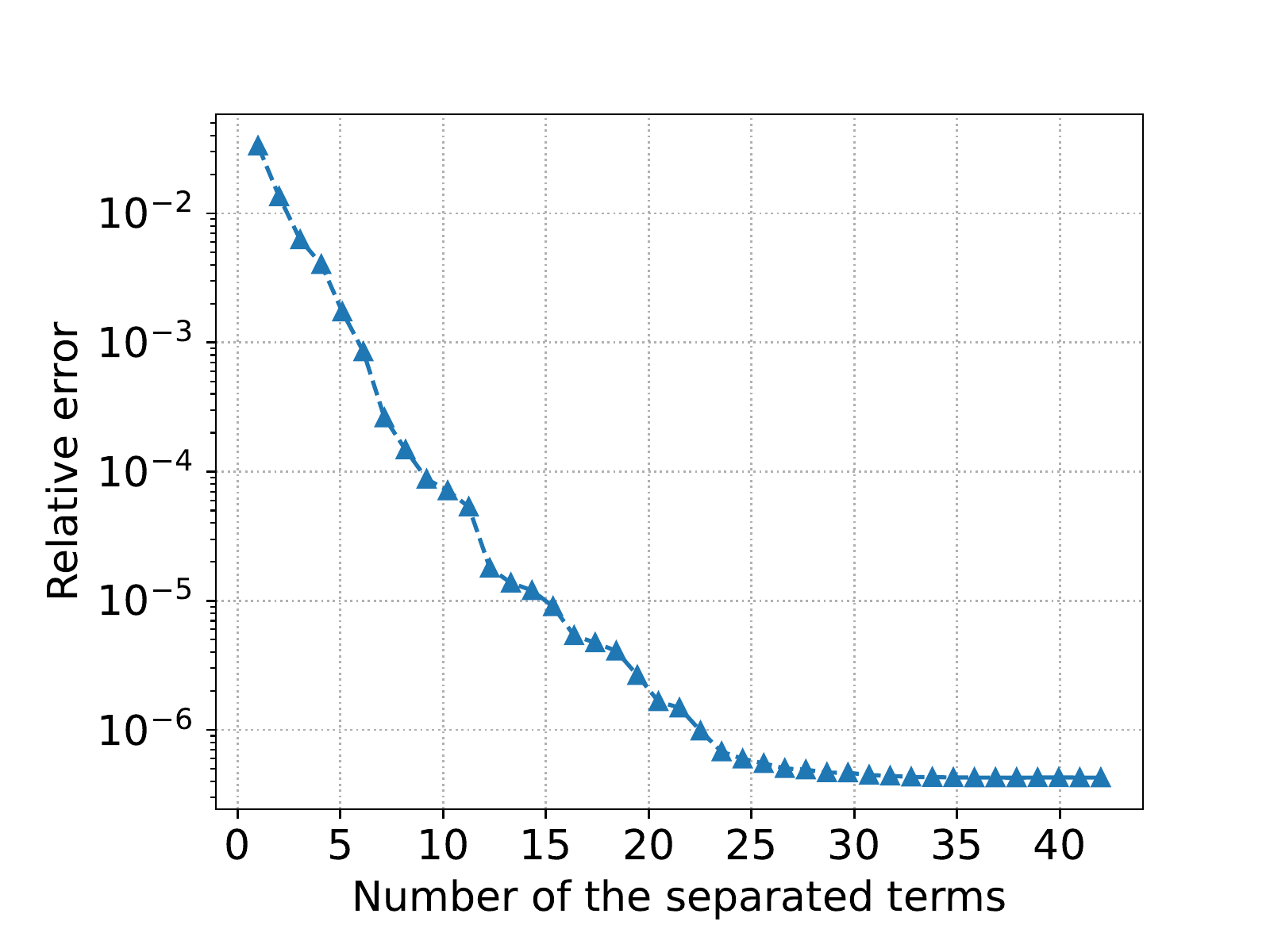}
	\caption{Comparison of the relative mean error corresponding to the different numbers of the separated terms $M$.}
	\label{fig3.1}
\end{figure}

First, we randomly choose $10^3$ samples and plot the relative mean error calculated by equation (\ref{eq-meanerror}) versus the number of the separated terms $M$ for the interface numerical solution in Figure \ref{fig3.1}, which shows that the SDD-VS method is suited for the high-dimensional stochastic problem, and the relative mean error becomes smaller when the number of the separated terms $M$ increases.

\begin{figure}[htbp]
	\centering
	\includegraphics[width=3.8in, height=2.7in]{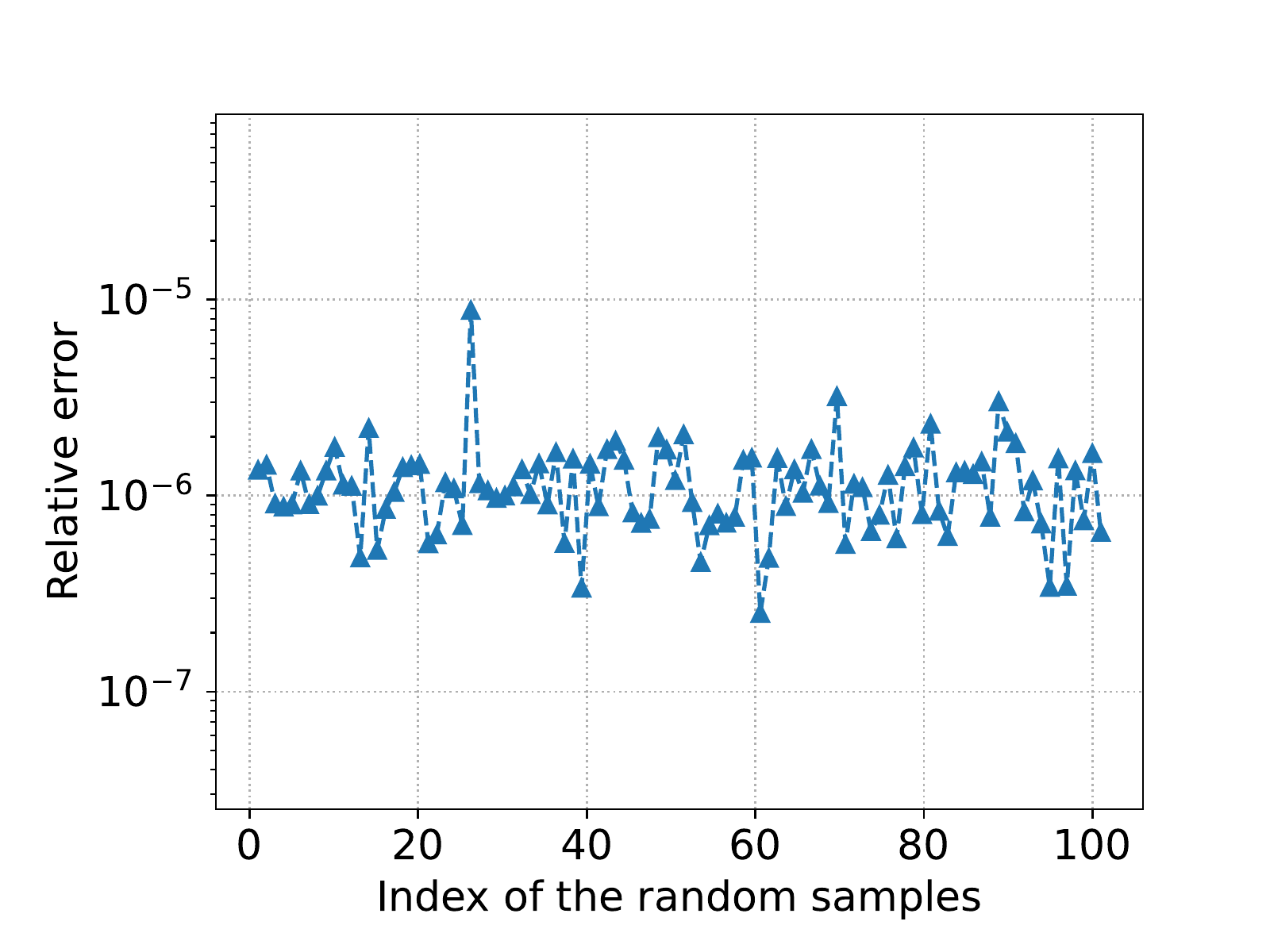}
	\caption{The relative mean error for 100 samples with the number of the separated terms $M=20$.}
	\label{fig3.2}
\end{figure}

We plot the relative error for the number of the separated terms $M=20$ in Figure \ref{fig3.2} to visualize the individual relative error by choosing the first $100$ samples out of $10^3$ random samples. From the figure, we can see that the SDD-VS method can give a good approximation for each sample.

\begin{figure}[htbp]
	\centering
	\includegraphics[width=3.8in, height=2.7in]{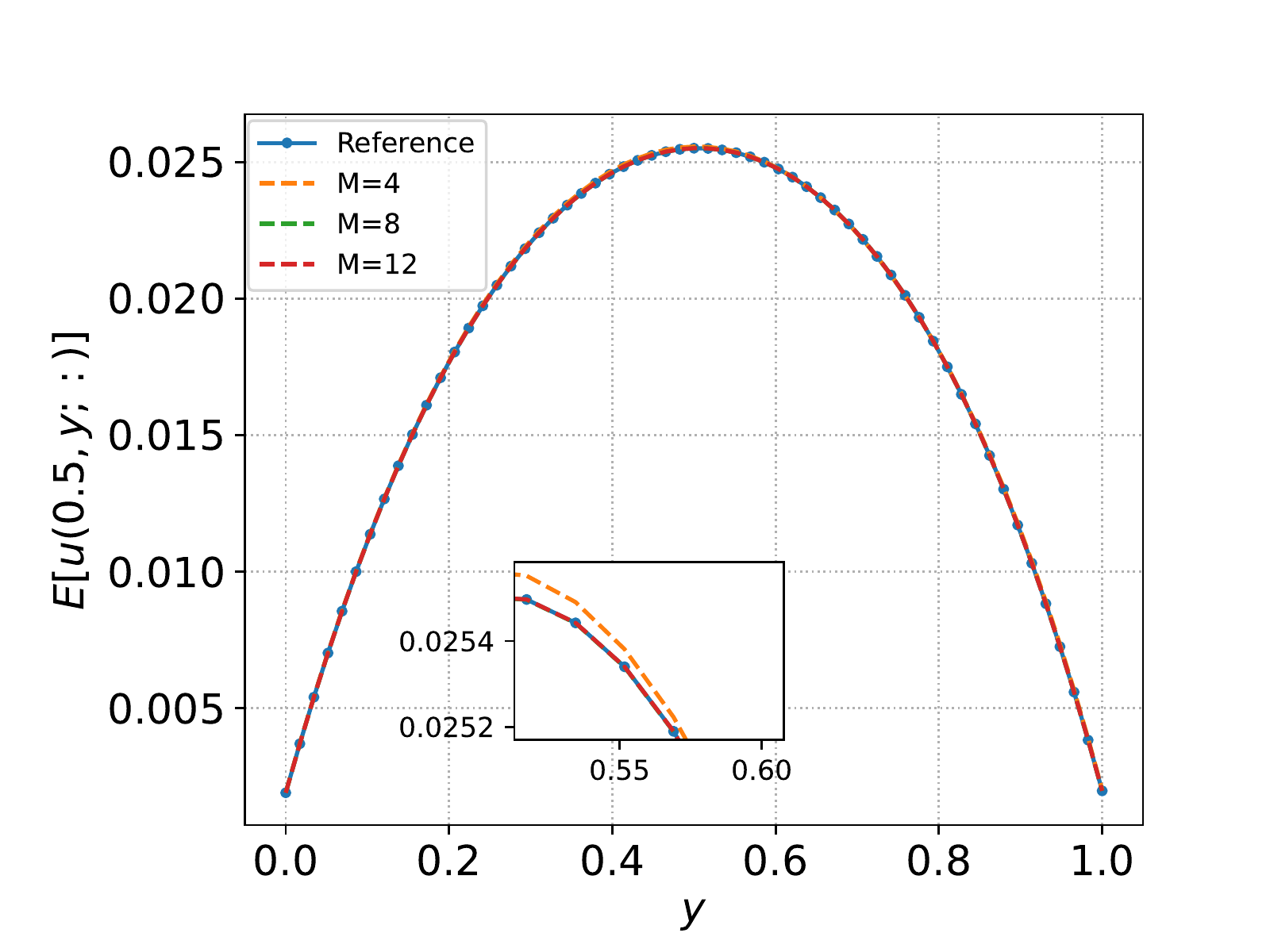}
	\caption{Comparison of the mean solution in the interface for reference and the SDD-VS method with the number of the separated terms $M=4,8,12$.}
	\label{fig3.3}
\end{figure}

In Figure \ref{fig3.3}, we plot the mean of the solution generated by our proposed SDD-VS method with the number of the separated terms $M=4,8,12$ and the reference in the interface based on $10^3$ random samples. As we can see, the approximation becomes better when the number of the separated terms $M$ increases.

\section{Conclusions}
\label{sec-Conclusions}
This paper presented the SDD-VS method for solving linear steady-state convection-diffusion equations with random coefficients.
The proposed method is devoted to building a relationship between the random inputs and the interface problem of SPDEs, in which the VS method plays a key role.
The whole computation of the SDD-VS method admits an offline–online decomposition.
In the offline phase,
we applied the direct Schur complement method to construct equation (\ref{eq-interface}), which is dense and depends on the inversion of $A_{II}^i(\bm{\xi})$, for $ 1\leq i \leq N_s$.
It brings great challenges to solving the stochastic interface system (\ref{eq-interface}) repeatedly, especially for many samples.
To improve computational efficiency, we adopted the extended VS method to reconstruct $S(\bm{\xi})$ and $F(\bm{\xi})$ as in equations (\ref{eq-S-F-Aff1}) and (\ref{eq-S-F-Aff2}), i.e., having affine decomposition. This renders the reduced model (\ref{eq-interf}), which needs much less computation effort than the original stochastic Schur complement system (\ref{eq-interface}) with the aid of affine decomposition. However, it may be not a very small-scale problem especially when the number of subdomains $N_s$ is large. To further improve efficiency, the VS method was used again to get the surrogate model for equation (\ref{eq-interf}).

In the online phase, we used the functional decomposition expression (\ref{eq-interu1})
to recover the solution to the stochastic interface problem for a large number of new samples. Moreover, the efficient surrogate model of the stochastic subproblem (\ref{eq-Sub-SPDE}) can be obtained by the VS method in \cite{Li2017a}. The online phase is efficient thanks to the reduced model representation for the stochastic interface solution. We applied the
proposed method to a few numerical models with random inputs. Careful numerical analysis was carried out for these numerical examples. We found that the SDD-VS method renders an efficient and robust reduced model.

In the future, we plan to extend the proposed method to tackle nonlinear unsteady problems. This will involve adapting the SDD-VS method to handle nonlinearity and time-dependent behavior. We also aim to conduct a rigorous convergence analysis for both the VS method and the SDD-VS method. This analysis will provide insights into the convergence properties and accuracy of the methods, enabling us to establish theoretical guarantees. Additionally, we will focus on applying the proposed method to models in dynamical systems, such as biological systems and petroleum engineering problems. By tailoring the methodology to these specific domains, we can address important challenges and develop efficient computational tools. In summary, our future research will involve extending the method to nonlinear unsteady problems, conducting convergence analysis, and exploring applications in dynamical systems. These efforts will contribute to advancing reduced-order modeling for stochastic partial differential equations and its practical use in various fields.

\end{document}